\documentclass[11pt]{article}
\usepackage{latexsym}

\newcommand{\eproof}{\mbox{\ }\hfill $\Box$ \par \vskip 10pt}

\newtheorem{Theorem}{Theorem}[section]
\newtheorem{lemma}[Theorem]{Lemma}
\newtheorem{prop}[Theorem]{Proposition}

\begin{document}

\title{Optimal dispersive estimates for the wave equation with $C^{\frac{n-3}{2}}$ potentials in dimensions $4\le n\le 7$}

\author{{\sc Fernando Cardoso and Georgi Vodev\thanks{Corresponding author}}}

\date{}

\maketitle

\noindent
{\bf Abstract.} We prove optimal dispersive estimates for the wave group $e^{it\sqrt{-\Delta+V}}$ for a class of real-valued potentials $V\in C^{\frac{n-3}{2}}({\bf R}^n)$, $4\le n\le 7$, such that $\partial_x^\alpha V(x)=O(\langle x\rangle^{-\delta})$, $\delta>\frac{n+1}{2}$, $|\alpha|\le \frac{n-3}{2}$. 

\setcounter{section}{0}
\section{Introduction and statement of results}

It is well known that the free wave group $e^{it\sqrt{G_0}}$ ($G_0$ being the self-adjoint realization of $-\Delta$ on $L^2({\bf R}^n)$, $n\ge 2$) satisfies the following dispersive estimates
$$\left\|e^{it\sqrt{G_0}}(\sqrt{G_0})^{-\frac{n+1}{2}+\epsilon}\langle G_0\rangle^{-\epsilon}\right\|_{L^1\to L^\infty}\le C_\epsilon|t|^{-\frac{n-1}{2}},\quad t\neq 0,
\eqno{(1.1)}$$
$$\left\|e^{it\sqrt{G_0}}(\sqrt{G_0})^{-\frac{n+1}{2}}\langle G_0\rangle^{-\epsilon}\right\|_{L^1\to L^\infty}\le C_\epsilon|t|^{-\frac{n-1}{2}}\log(2+|t|),\quad t\neq 0,
\eqno{(1.2)}$$
for every $0<\epsilon\ll 1$, and
$$\left\|e^{it\sqrt{G_0}}(\sqrt{G_0})^{-\frac{\alpha(n+1)}{2}}\right\|_{L^{p'}\to L^p}\le C|t|^{-\frac{\alpha(n-1)}{2}},\quad t\neq 0,
\eqno{(1.3)}$$
for every $2\le p<+\infty$, where $1/p+1/p'=1$ and $\alpha=1-2/p$. Note that (1.1) and (1.2) do not
hold with $\epsilon=0$.

The problem we address in the present paper is that one of finding as large as possible class of real-valued potentials, $V$, such that the self-adjoint realization, $G$, of the
operator $-\Delta+V$ on $L^2({\bf R}^n)$ satisfies estimates similar to (1.1)-(1.3). In dimensions two and three this problem is actually solved and in particular one knows that no regularity of the potential is needed in order to have analogues of (1.1)-(1.3) for the operator $G$ (see \cite{kn:CCV1}, \cite{kn:M2}, \cite{kn:DP}, \cite{kn:GV}). The same conlusion remains true 
in higher dimensions as far as the low and the intermediate frequencies are concerned (see \cite{kn:M}, \cite{kn:V}), while at high frequencies one is obliged to loose derivatives if no
regularity of the potential is required. Indeed, dispersive estimates with a loss of $\frac{n-3}{2}$ derivatives for the perturbed wave group were proved in \cite{kn:V} for potentials $V\in L^\infty({\bf R}^n)$ satisfying
$$\left|V(x)\right|\le C\langle x\rangle^{-\delta},\quad\forall x\in{\bf R}^n,\eqno{(1.4)}$$
with constants $C>0$ and $\delta>\frac{n+1}{2}$. In other words, to get optimal dispersive estimates for the perturbed wave group when $n\ge 4$ one needs to assume some regularity on the potential. 
Indeed, such estimates were proved in \cite{kn:B} for potentials belonging to the Schwartz class. 
Getting the minimal regularity of the potential in order to have optimal dispersive estimates for the perturbed wave group when $n\ge 4$, however, turns out to be a hard problem. The counterexample of \cite{kn:GoV} shows the existence of potentials $V\in C_0^k({\bf R}^n)$, $\forall k<\frac{n-3}{2}$, for which the perturbed Schr\"odinger group $e^{itG}$ does not
satisfy optimal $L^1\to L^\infty$ dispersive estimates. In analogy, one could expect that a
similar phenomenon occurs for the wave group, too. Thus the natural question is to ask if we have 
optimal dispersive estimates for potentials $V\in C^{\frac{n-3}{2}}({\bf R}^n)$, $n\ge 4$.
Indeed, such a result has been recently proved in \cite{kn:EG} for the Schr\"odinger group $e^{itG}$ when $n=5,7$, while in \cite{kn:CCV2} this was previously proved for potentials $V\in C^k({\bf R}^n)$, $ k>\frac{n-3}{2}$, $n=4,5$. Let us also mention the work \cite{kn:CCV3} where
$L^1\to L^\infty$ dispersive estimates for $e^{itG}$ with a logarithmic loss of derivatives were proved for potentials $V\in C^{\frac{n-3}{2}}({\bf R}^n)$ still in dimensions four and five.  
To our best knowledge, no such results exist for the perturbed wave group $e^{it\sqrt{G}}$. The purpose of this work is to prove this when $4\le n\le 7$. To be more precise, 
define the sets of functions ${\cal V}^k_\delta({\bf R}^n)$, ${\cal C}^k_\delta({\bf R}^n)$, $\delta,k\ge 0$, as follows. If $k$ is integer, ${\cal V}^k_\delta({\bf R}^n)$ (res. ${\cal C}^k_\delta({\bf R}^n)$) is the set of all functions $V\in C^k({\bf R}^n)$ satisfying respectively
$$\left\|V\right\|_{{\cal V}_\delta^k}:=\sup_{x\in{\bf R}^n}\sum_{0\le|\alpha|\le k}\langle x\rangle^{\delta+|\alpha|}\left|\partial_x^\alpha V(x)\right|<+\infty,$$ 
$$\left\|V\right\|_{{\cal C}_\delta^k}:=\sup_{x\in{\bf R}^n}\sum_{0\le|\alpha|\le k}\langle x\rangle^{\delta}\left|\partial_x^\alpha V(x)\right|<+\infty.$$ 
If $k=k_0+\nu$ with $k_0\ge 0$ an integer and $0<\nu<1$, a function $V$ will be said to belong to ${\cal V}^k_\delta({\bf R}^n)$ (res. ${\cal C}^k_\delta({\bf R}^n)$) if $V\in{\cal V}^{k_0}_\delta({\bf R}^n)$ (res. ${\cal C}^{k_0}_\delta({\bf R}^n)$) and if there exists a family of functions 
$V_\theta\in{\cal V}^{k_0+1}_\delta({\bf R}^n)$ (res. ${\cal C}^{k_0+1}_\delta({\bf R}^n)$), $0<\theta\le 1$, such that
$$\left\|V\right\|_{{\cal V}_\delta^k({\cal C}_\delta^k)}:=\left\|V\right\|_{{\cal V}_\delta^{k_0}({\cal C}_\delta^{k_0})}+\sup_{0<\theta\le 1}\left(\theta^{-\nu}\left\|V-V_\theta\right\|_{{\cal V}_\delta^{k_0}({\cal C}_\delta^{k_0})} +
\theta^{1-\nu}\left\|V_\theta\right\|_{{\cal V}_\delta^{k_0+1}({\cal C}_\delta^{k_0+1})}\right)<+\infty.$$
Our conjecture is that the perturbed wave group $e^{it\sqrt{G}}$ satisfies optimal dispersive estimates for real-valued potentials satisfying (1.4) as well as the condition
$$V\in {\cal V}^{\frac{n-3}{2}}_\delta({\bf R}^n),\quad \delta>2,\eqno{(1.5)}$$
(or probably only (1.5)).
In the present paper we prove optimal dispersive estimates when $4\le n\le 7$ under the following stronger condition:
$$V\in {\cal C}^{\frac{n-3}{2}}_\delta({\bf R}^n),\quad \delta>\frac{n+1}{2}.\eqno{(1.6)}$$
Proving this when $n\ge 8$, however, remains an open problem. 
Given $a>0$, choose a function
$\chi_a\in C^\infty({\bf R})$, $\chi_a(\lambda)=0$ for $\lambda\le a$, $\chi_a(\lambda)=1$ for $\lambda\ge a+1$. Our main result is the following

\begin{Theorem} Let $4\le n\le 7$ and suppose that $V$ satisfies (1.6). Then, for every 
$a>0$, $0<\epsilon\ll 1$, $2\le p<+\infty$, $t\neq 0$, we have the estimates
$$\left\|e^{it\sqrt{G}}(\sqrt{G})^{-\frac{n+1}{2}-\epsilon}\chi_a(\sqrt{G})\right\|_{L^1\to L^\infty}\le C_{\epsilon,a}|t|^{-\frac{n-1}{2}},\eqno{(1.7)}$$
$$\left\|e^{it\sqrt{G}}(\sqrt{G})^{-\frac{\alpha(n+1)}{2}}\chi_a(\sqrt{G})\right\|_{L^{p'}\to L^p}\le C_a|t|^{-\frac{\alpha(n-1)}{2}},\eqno{(1.8)}$$
 where $1/p+1/p'=1$ and $\alpha=1-2/p$. Moreover, if in addition we suppose that zero is neither
 an eigenvalue nor a resonance of $G$, then we have the estimates
 $$\left\|e^{it\sqrt{G}}(\sqrt{G})^{-\frac{n+1}{2}+\epsilon}\langle G\rangle^{-\epsilon}P_{ac}\right\|_{L^1\to L^\infty}\le C_{\epsilon}|t|^{-\frac{n-1}{2}},\eqno{(1.9)}$$
 $$\left\|e^{it\sqrt{G}}(\sqrt{G})^{-\frac{n+1}{2}}\langle G\rangle^{-\epsilon}P_{ac}\right\|_{L^1\to L^\infty}\le C_{\epsilon}|t|^{-\frac{n-1}{2}}\log(2+|t|),\eqno{(1.10)}$$
$$\left\|e^{it\sqrt{G}}(\sqrt{G})^{-\frac{\alpha(n+1)}{2}}P_{ac}\right\|_{L^{p'}\to L^p}\le C|t|^{-\frac{\alpha(n-1)}{2}},\eqno{(1.11)}$$
where $P_{ac}$ denotes the spectral projection onto the absolutely continuous spectrum of $G$.
\end{Theorem}

\noindent
{\bf Remark.} In view of the low frequency dispersive estimates proved in \cite{kn:M} under
the assumption (1.4), the estimates (1.9) and (1.10) follow from (1.7), while (1.11) follows from (1.8).
It is worth also noticing that it suffices to prove (1.7) and (1.8) for $a\gg 1$ as at intermediate frequencies the dispersive estimates are proved in \cite{kn:V} under (1.4) only.

To prove (1.7) and (1.8) it suffices to prove an almost optimal (in $h$) bound of the $L^1\to L^\infty$ norm of the operator  $e^{it\sqrt{G}}\varphi(h\sqrt{G})$, where
$\varphi\in C_0^\infty((0,+\infty))$ and $0<h\ll 1$ (see Theorem 2.1 below). 
Then we reduce this problem (using only (1.4)) to estimating the $L^1\to L^\infty$ norm of $\left[\frac{n-2}{2}\right]$ operators (denoted by ${\cal A}_k$ below) with explicit kernels (see Theorem 2.2 below). In particular, when $n=4,5$ one needs to estimate the $L^1\to L^\infty$ norm of only one operator, ${\cal A}_1$ (see Section 4), while
when $n=6,7$ one must also bound the $L^1\to L^\infty$ norm of another operator, ${\cal A}_2$ (see Section 5). In higher dimensions the kernels get much more complicated and therefore the problem gets much more technical and harder. In fact, the kernels of ${\cal A}_k$ are oscilatory integrals with $O\left(h^{-1}\right)$ non-smooth phases. Thus, the only way to gain behavior in $h$ is the integration by parts with respect to the space variables - roughly speaking, one needs to integrate by parts $\frac{(n-3)k}{2}$ times the kernel of ${\cal A}_k$. This, however, leads to singular integrals, so the most delicate point of the proof consists of finding such an integration by parts scheme that allows to avoid non-integrable singularities. When $k\ge 2$ this turns out to be a very tough problem with an increasing complexity as $n$ grows up, as indicated in \cite{kn:EG} in the context of the Schr\"odinger equation. In contrast, the analysis of the operator 
${\cal A}_1$ is relatively easy and can be carried out in all dimensions $n\ge 4$. Note also that a more complicated integration by parts scheme than that one we use (see the proof of Proposition 5.5) could probably lead to the estimate (2.1) below with $\epsilon=0$ (when $n=7$). However, since the $\epsilon$ loss in (2.1) does not affect the proof of the main result, we prefer to keep the proof relatively simple and short allowing an $\epsilon$ loss in $h$ (in the case $n=7$) rather than seeking a sharp estimate by much more complicated arguments. 
We finally reduce the problem to bounding singular integrals essentially studied in Section 6 of \cite{kn:EG}. The bounds we need are actually simpler than those proved in \cite{kn:EG} - we sketch the proof for the sake of completeness in the appendix of the present paper.  
Note that our method works also in even dimensions, though the proof is more difficult. Indeed, in this case our estimates can be proved by applying interpolation arguments (relatively easy when $n=4$ and much more complicated when $n=6$).

\section{Reduction to semi-classical dispersive estimates}

We will first show that the estimates (1.7) and (1.8) follow from the following

\begin{Theorem} Let $\varphi\in C_0^\infty((0,+\infty))$. Then, under the assumptions of Theorem 1.1, for all $0<h\le 1$, $t\neq 0$, $0<\epsilon\ll 1$, we have the estimate
$$\left\|e^{it\sqrt{G}}\varphi(h\sqrt{G})\right\|_{L^1\to L^\infty}\le C_\epsilon h^{-\frac{n+1}{2}-\epsilon}|t|^{-\frac{n-1}{2}},\eqno{(2.1)}$$
with a constant $C_\epsilon>0$ independent of $h$ and $t$.
\end{Theorem}

\noindent
{\bf Remark.} When $4\le n\le 6$, we actually prove (2.1) with $\epsilon=0$.

To prove (1.7) we will use the identity
$$\sigma^{-\frac{n+1}{2}-\epsilon}\chi_a(\sigma)=\int_0^1\varphi(h\sigma)h^{\frac{n+1}{2}+\epsilon-1}dh,$$
where $\varphi(\sigma)=\sigma^{1-\frac{n+1}{2}-\epsilon}\chi'_a(\sigma)\in C_0^\infty((0,+\infty))$. Using (2.1) we get
$$\left\|e^{it\sqrt{G}}(\sqrt{G})^{-\frac{n+1}{2}-\epsilon}\chi_a(\sqrt{G})\right\|_{L^1\to L^\infty}\le \int_0^1\left\|e^{it\sqrt{G}}\varphi(h\sqrt{G})\right\|_{L^1\to L^\infty}h^{\frac{n+1}{2}+\epsilon-1}dh$$ $$\le C_\epsilon|t|^{-\frac{n-1}{2}}\int_0^1h^{-1+\frac{\epsilon}{2}}dh\le
C_\epsilon |t|^{-\frac{n-1}{2}}.$$
To prove (1.8) we will use the identity
$$\sigma^{-\frac{\alpha(n+1)}{2}}\chi_a(\sigma)=\int_0^1\varphi(h\sigma)h^{\frac{\alpha(n+1)}{2}-1}dh,$$
where $\varphi(\sigma)=\sigma^{1-\frac{\alpha(n+1)}{2}}\chi'_a(\sigma)\in C_0^\infty((0,+\infty))$. Since the operator $G_0$ satisfies (2.1), by Theorem 2.1 we get
$$\left\|e^{it\sqrt{G}}\varphi(h\sqrt{G})-e^{it\sqrt{G_0}}\varphi(h\sqrt{G_0})\right\|_{L^1\to L^\infty}\le C_\epsilon h^{-\frac{n+1}{2}-\epsilon}|t|^{-\frac{n-1}{2}}.\eqno{(2.2)}$$
On the other hand, we have (see Theorem 3.1 of \cite{kn:V})
$$\left\|e^{it\sqrt{G}}\varphi(h\sqrt{G})-e^{it\sqrt{G_0}}\varphi(h\sqrt{G_0})\right\|_{L^2\to L^2}\le Ch,\quad\forall\,t.\eqno{(2.3)}$$
By interpolation between (2.2) and (2.3) we conclude
$$\left\|e^{it\sqrt{G}}\varphi(h\sqrt{G})-e^{it\sqrt{G_0}}\varphi(h\sqrt{G_0})\right\|_{L^{p'}\to L^p}\le C_\epsilon h^{1-\frac{\alpha(n+3)}{2}-\epsilon\alpha}|t|^{-\frac{\alpha(n-1)}{2}},\eqno{(2.4)}$$
for every $2\le p\le +\infty$, where $1/p+1/p'=1$, $\alpha=1-2/p$. Using (2.4) we get
 $$\left\|e^{it\sqrt{G}}(\sqrt{G})^{-\frac{\alpha(n+1)}{2}}\chi_a(\sqrt{G})-e^{it\sqrt{G_0}}
 (\sqrt{G_0})^{-\frac{\alpha(n+1)}{2}}\chi_a(\sqrt{G_0})\right\|_{L^{p'}\to L^p}$$
 $$\le\int_0^1\left\|e^{it\sqrt{G}}\varphi(h\sqrt{G})-e^{it\sqrt{G_0}}\varphi(h\sqrt{G_0})\right\|_{L^{p'}\to L^p}h^{\frac{\alpha(n+1)}{2}-1}dh$$
 $$\le C_\epsilon|t|^{-\frac{\alpha(n-1)}{2}}\int_0^1h^{-\alpha(1+\epsilon)}dh\le
C'|t|^{-\frac{\alpha(n-1)}{2}},\eqno{(2.5)}$$
provided $2\le p<+\infty$ (i.e. $0\le\alpha<1$) and $\epsilon$ is such that $\alpha(1+\epsilon)<1$. Clearly, (1.8) follows from (2.5) and (1.3).\\

{\it Proof of Theorem 2.1.} Note that it suffices to prove (2.1) for $0<h\le h_0$ with some constant $0<h_0\ll 1$ as for $h_0\le h\le 1$
it is proved in \cite{kn:V} under (1.4) only. We are going to use the formula
$$e^{it\sqrt{G_0}}\varphi(h\sqrt{G_0})=(\pi i)^{-1}\int_0^\infty e^{it\lambda}\varphi(h\lambda)\left(R_0^+(\lambda)-R_0^-(\lambda)\right)\lambda
d\lambda,\eqno{(2.6)}$$
where $R_0^\pm(\lambda)=(G_0-\lambda^2\pm i0)^{-1}$ are the outgoing and incoming free resolvents with kernels given by
$$[R_0^\pm(\lambda)](x,y)=\frac{\pm i|x-y|^{-2\nu}}{4(2\pi)^\nu}{\cal H}_\nu^\pm(\lambda|x-y|),$$
where $\nu=\frac{n-2}{2}$, ${\cal H}_\nu^\pm(z)=z^\nu H_\nu^\pm(z)$, $H_\nu^\pm(z)$ being the Hankel functions of order $\nu$. We also have the formula
$$e^{it\sqrt{G}}\varphi(h\sqrt{G})=(\pi i)^{-1}\int_0^\infty e^{it\lambda}\varphi(h\lambda)\left(R^+(\lambda)-R^-(\lambda)\right)\lambda
d\lambda,\eqno{(2.7)}$$
where $R^\pm(\lambda)=(G-\lambda^2\pm i0)^{-1}$ are the outgoing and incoming perturbed resolvents satisfying the relation
$$R^\pm(\lambda)\left(1+VR_0^\pm(\lambda)\right)=R_0^\pm(\lambda).\eqno{(2.8)}$$
Iterating (2.8) $m$ times we get the identity
$$R^\pm(\lambda)-R_0^\pm(\lambda)=\sum_{k=1}^mR_0^\pm(\lambda)\left(-VR_0^\pm(\lambda)\right)^k+R^\pm(\lambda)\left(-VR_0^\pm(\lambda)\right)^{m+1}.\eqno{(2.9)}$$
In view of (2.6), (2.7) and (2.9), we can write
$$e^{it\sqrt{G}}\varphi(h\sqrt{G})-e^{it\sqrt{G_0}}\varphi(h\sqrt{G_0})=\sum_{k=1}^m{\cal A}_k(t,h)+{\cal R}_m(t,h),\eqno{(2.10)}$$
where ${\cal A}_k={\cal A}_k^+-{\cal A}_k^-$, ${\cal R}_m={\cal R}_m^+-{\cal R}_m^-$,
$${\cal A}_k^\pm(t,h)=(h\pi i)^{-1}\int_0^\infty e^{it\lambda}\widetilde\varphi(h\lambda)R_0^\pm(\lambda)\left(-VR_0^\pm(\lambda)\right)^k
d\lambda,$$
$${\cal R}_m^\pm(t,h)=(h\pi i)^{-1}\int_0^\infty e^{it\lambda}\widetilde\varphi(h\lambda)R^\pm(\lambda)\left(-VR_0^\pm(\lambda)\right)^{m+1}
d\lambda,$$
where $\widetilde\varphi(\lambda)=\lambda\varphi(\lambda)$. In the next section we will prove the following

\begin{Theorem} Suppose that $V$ satisfies (1.4). Then, in all dimensions $n\ge 4$ and for all $0<h\ll 1$,
$t\neq 0$, we have the estimates
$$\left\|e^{it\sqrt{G}}\varphi(h\sqrt{G})-e^{it\sqrt{G_0}}\varphi(h\sqrt{G_0})-\sum_{k=1}^{\left[\frac{n-2}{2}\right]}{\cal A}_k(t,h)\right\|_{L^1\to L^\infty}\le Ch^{-\frac{n+1}{2}}|t|^{-\frac{n-1}{2}},\eqno{(2.11)}$$
$$\left\|{\cal A}_k(t,h)\right\|_{L^1\to L^\infty}\le C_kh^{k-n}|t|^{-\frac{n-1}{2}},\quad \forall\,k\ge 1.\eqno{(2.12)}$$
\end{Theorem}
 
Thus, to prove (2.1) when $n=4,5$ it suffices to improve (2.12) in $h$ for $k=1$, only, using the regularity assumption (1.5), while when $n=6,7$ it suffices to improve (2.12) for $k=1,2$, using (1.6). This analysis will be carried out in Sections 4 and 5. Note also that in dimensions $n=2,3$ the above theorem is proved in \cite{kn:CCV1}.

\section{Proof of Theorem 2.2}
Set
$$T_k^\pm(\lambda)=R_0^\pm(\lambda)\left(-VR_0^\pm(\lambda)\right)^k,\quad
\widetilde T_k^\pm(\lambda)=R^\pm(\lambda)\left(-VR_0^\pm(\lambda)\right)^k.$$
We will first show that Theorem 2.2 follows from the following

\begin{prop} Under the assumptions of Theorem 2.2, there exists a constant $\lambda_0>0$ so that if $n$ is odd, for all integers $k\ge 1$, $0\le m\le\frac{n-1}{2}$, we have
$$\left\|\frac{d^mT_k^\pm(\lambda)}{d\lambda^m}\right\|_{L^1\to L^\infty}\le C_k\lambda^{n-2-k},\eqno{(3.1)}$$
for $\lambda\ge\lambda_0$. If $n$ is even, (3.1) still holds for $0\le m\le\frac{n-2}{2}$. Moreover, in this case we also have
$$\left\|\frac{d^{\frac{n-2}{2}}T_k^\pm}{d\lambda^{\frac{n-2}{2}}}(\lambda_1)-\frac{d^{\frac{n-2}{2}}T_k^\pm}{d\lambda^{\frac{n-2}{2}}}(\lambda_2)\right\|_{L^1\to L^\infty}\le C_k\lambda_1^{n-2-k}|\lambda_1-\lambda_2|^{1/2},\eqno{(3.2)}$$
for $\lambda_1+1\ge\lambda_2>\lambda_1\ge\lambda_0$. These estimates remain valid with $T_k^\pm$ replaced by $\widetilde T_k^\pm$.
\end{prop}

If $n$ is odd, then $\left[\frac{n-2}{2}\right]=\frac{n-3}{2}$. In this case (2.12) follows from (3.1) by integrating by parts $\frac{n-1}{2}$
times with respect to the variable $\lambda$. Similarly, (2.11) follows from (2.10) applied with $m=\frac{n-3}{2}$ and (3.1) used with $T_k^\pm$ replaced by $\widetilde T_k^\pm$, $k=\frac{n-1}{2}$. Let now $n$ be even. Choose a real-valued function $\widetilde\phi\in C_0^\infty([1,2])$, $\widetilde\phi\ge 0$, such that $\int\widetilde\phi(\sigma)d\sigma=1$ and set
$$T_{k,\theta}^\pm(\lambda)=\theta^{-1}\int T_k^\pm(\lambda+\sigma)\widetilde\phi(\sigma/\theta)d\sigma,\quad 0<\theta\le 1.$$
It follows from Proposition 3.1 that for $0\le m\le\frac{n-2}{2}$, $\lambda\ge\lambda_0$, we have
$$\left\|\frac{d^mT_{k,\theta}^\pm(\lambda)}{d\lambda^m}\right\|_{L^1\to L^\infty}\le C_k\lambda^{n-2-k},\eqno{(3.3)}$$
$$\left\|\frac{d^m(T_{k,\theta}^\pm-T_k^\pm)(\lambda)}{d\lambda^m}\right\|_{L^1\to L^\infty}\le C_k\theta^{1/2}\lambda^{n-2-k},\eqno{(3.4)}$$
$$\left\|\frac{d^\frac{n}{2}T_{k,\theta}^\pm(\lambda)}{d\lambda^\frac{n}{2}}\right\|_{L^1\to L^\infty}\le  C_k\theta^{-1/2}\lambda^{n-2-k}.\eqno{(3.5)}$$
Integrating by parts $\frac{n-2}{2}$ times and using (3.4), we get
$$\left\|(h\pi i)^{-1}\int_0^\infty e^{it\lambda}\widetilde\varphi(h\lambda)\left(T_{k,\theta}^\pm(\lambda)-T_k^\pm(\lambda)\right)d\lambda\right\|_{L^1\to L^\infty}\le C_k\theta^{1/2}h^{k-n}|t|^{-\frac{n-2}{2}}.\eqno{(3.6)}$$
Integrating by parts $\frac{n}{2}$ times and using (3.3) and (3.5), we get
$$\left\|(h\pi i)^{-1}\int_0^\infty e^{it\lambda}\widetilde\varphi(h\lambda)T_{k,\theta}^\pm(\lambda)d\lambda\right\|_{L^1\to L^\infty}\le C_k\theta^{-1/2}h^{k-n}|t|^{-\frac{n}{2}}.\eqno{(3.7)}$$
By (3.6) and (3.7),
$$\left\|(h\pi i)^{-1}\int_0^\infty e^{it\lambda}\widetilde\varphi(h\lambda)T_{k}^\pm(\lambda)d\lambda\right\|_{L^1\to L^\infty}\le C_kh^{k-n}|t|^{-\frac{n-1}{2}}\left((\theta|t|)^{1/2}+(\theta|t|)^{-1/2}\right).\eqno{(3.8)}$$
If $|t|\ge 1$ we take $\theta=|t|^{-1}$ in (3.8) to conclude
$$\left\|{\cal A}_k^\pm(t,h)\right\|_{L^1\to L^\infty}\le C_kh^{k-n}|t|^{-\frac{n-1}{2}}.\eqno{(3.9)}$$
If $|t|\le 1$ the estimate (3.9) follows directly from (3.1) with $m=0$ without integrating by parts. The estimate (2.11) follows in precisely
the same way using (2.10) with $m=\frac{n-2}{2}$ and replacing $T_k^\pm$ by $\widetilde T_k^\pm$.\\

{\it Proof of Proposition 3.1.} It is well known that the functions ${\cal H}^\pm_\nu$ satisfy
$$\partial_z^k{\cal H}_\nu^\pm(z)=O_k\left(z^{\frac{n-3}{2}}\right),\quad z\ge 1,\eqno{(3.10)}$$ for all integers $k\ge 0$, while at $z=0$ they 
are of the form
 $${\cal H}^\pm_\nu(z)={\cal H}^\pm_{\nu,1}(z)+z^{n-2}\log z{\cal H}^\pm_{\nu,2}(z),\eqno{(3.11)}$$
 where ${\cal H}^\pm_{\nu,j}$ are analytic, ${\cal H}^\pm_{\nu,2}\equiv 0$ if $n$ is odd.
Let $\phi\in C_0^\infty({\bf R})$, $\phi(z)=1$ for $|z|\le 1/2$, $\phi(z)=0$ for $|z|\ge 1$. Decompose the operator $R_0^\pm(\lambda)$ as ${\cal K}_1^\pm(\lambda)+{\cal K}_2^\pm(\lambda)$, where the kernels of ${\cal K}_1^\pm(\lambda)$ and ${\cal K}_2^\pm(\lambda)$ are defined by replacing in the kernel of $R_0^\pm(\lambda)$ the function ${\cal H}_\nu^\pm$ by $(1-\phi){\cal H}_\nu^\pm$ and $\phi{\cal H}_\nu^\pm$, respectively.

\begin{lemma} If $n$ is odd, for all integers $0\le m\le\frac{n-1}{2}$ and all $\lambda\ge 1$, $0<\epsilon\ll 1$, we have
$$\left\|\frac{d^m{\cal K}_1^\pm(\lambda)}{d\lambda^m}\langle x\rangle^{-1/2-m-\epsilon}\right\|_{L^2\to L^\infty}+
\left\|\langle x\rangle^{-1/2-m-\epsilon}\frac{d^m{\cal K}_1^\pm(\lambda)}{d\lambda^m}\right\|_{L^1\to L^2}\le C\lambda^{\frac{n-3}{2}},\eqno{(3.12)}$$
$$\left\|\frac{d^m{\cal K}_2^\pm(\lambda)}{d\lambda^m}\langle x\rangle^{-1-m-\epsilon}\right\|_{L^\infty\to L^\infty}+
\left\|\langle x\rangle^{-1-m-\epsilon}\frac{d^m{\cal K}_2^\pm(\lambda)}{d\lambda^m}\right\|_{L^1\to L^1}\le C\lambda^{-1}.\eqno{(3.13)}$$
If $n$ is even, (3.12) and (3.13) still hold for $0\le m\le\frac{n-2}{2}$. In this case we also have with $m=\frac{n-2}{2}$
$$\left\|\left(\frac{d^m{\cal K}_1^\pm}{d\lambda^m}(\lambda_1)-\frac{d^m{\cal K}_1^\pm}{d\lambda^m}(\lambda_2)\right)\langle x\rangle^{-1-m-\epsilon}\right\|_{L^2\to L^\infty}$$ $$+
\left\|\langle x\rangle^{-1-m-\epsilon}\left(\frac{d^m{\cal K}_1^\pm}{d\lambda^m}(\lambda_1)-\frac{d^m{\cal K}_1^\pm}{d\lambda^m}(\lambda_2)\right)\right\|_{L^1\to L^2}\le C\lambda_1^{\frac{n-3}{2}}|\lambda_1-\lambda_2|^{1/2},\eqno{(3.14)}$$
$$\left\|\left(\frac{d^m{\cal K}_2^\pm}{d\lambda^m}(\lambda_1)-\frac{d^m{\cal K}_2^\pm}{d\lambda^m}(\lambda_2)\right)\langle x\rangle^{-3/2-m-\epsilon}\right\|_{L^\infty\to L^\infty}$$ $$+
\left\|\langle x\rangle^{-3/2-m-\epsilon}\left(\frac{d^m{\cal K}_2^\pm}{d\lambda^m}(\lambda_1)-\frac{d^m{\cal K}_2^\pm}{d\lambda^m}(\lambda_2)\right)\right\|_{L^1\to L^1}\le C\lambda_1^{-1}|\lambda_1-\lambda_2|^{1/2},\eqno{(3.15)}$$
for $\lambda_1+1\ge\lambda_2>\lambda_1\ge 1$.
\end{lemma}

{\it Proof.} Denote by $K_{j,m}^\pm(x,y,\lambda)$ the kernel of the operator $\frac{d^m{\cal K}_j^\pm(\lambda)}{d\lambda^m}$, $j=1,2$. In view of (3.10), we have
$$\left|K_{1,m}^\pm(x,y,\lambda)\right|\le C\lambda^{\frac{n-3}{2}}|x-y|^{m-\frac{n-1}{2}}.\eqno{(3.16)}$$
On the other hand, it is easy to see that the left-hand side of (3.12) is equivalent to the square root of
$$\sup_{y\in{\bf R}^n}\int_{{\bf R}^n}\left|K_{1,m}^\pm(x,y,\lambda)\right|^2\langle x\rangle^{-1-2m-2\epsilon}dx$$ $$\le C'\lambda^{n-3}\sup_{y\in{\bf R}^n}\int_{{\bf R}^n}|x-y|^{2m-n+1}\langle x\rangle^{-1-2m-2\epsilon}dx\le C\lambda^{n-3},$$
provided $0\le m\le\left[\frac{n-1}{2}\right]$. For these values of $m$ we also have, in view of (3.11),
$$\left|K_{2,m}^\pm(x,y,\lambda)\right|\le C\lambda^{-1}|x-y|^{m-n+1}.\eqno{(3.17)}$$
Thus, the left-hand side of (3.13) is equivalent to  
$$\sup_{y\in{\bf R}^n}\int_{{\bf R}^n}\left|K_{2,m}^\pm(x,y,\lambda)\right|\langle x\rangle^{-1-m-\epsilon}dx$$ $$\le C'\lambda^{-1}\sup_{y\in{\bf R}^n}\int_{{\bf R}^n}|x-y|^{m-n+1}\langle x\rangle^{-1-m-\epsilon}dx\le C\lambda^{-1}.$$
To prove (3.14) we will use that given any function $f\in C^1({\bf R})$ and any $\sigma>0$ we have the inequality
$$|f(\sigma\lambda_1)-f(\sigma\lambda_2)|^2\le\sigma(|f(\sigma\lambda_1)|+|f(\sigma\lambda_2)|)\int_{\lambda_1}^{\lambda_2}|f'(\sigma\lambda)|d\lambda$$ $$\le\sigma|\lambda_1-\lambda_2|\left(|f(\sigma\lambda_1)|^2+|f(\sigma\lambda_2)|^2\right)+\sigma\int_{\lambda_1}^{\lambda_2}|f'(\sigma\lambda)|^2d\lambda,\eqno{(3.18)}$$
where $f'(z)=df(z)/dz$. Applying (3.18) with $f(z)=\frac{d^m((1-\phi(z)){\cal H}_\nu^\pm(z))}{dz^m}$, $m=\frac{n-2}{2}$, $\sigma=|x-y|$, and using  (3.10), we obtain
$$\left|K_{1,m}^\pm(x,y,\lambda_1)-K_{1,m}^\pm(x,y,\lambda_2)\right|^2\le C\lambda_1^{n-3}|\lambda_1-\lambda_2|.\eqno{(3.19)}$$
Hence, the left-hand side of (3.14) is equivalent to the square root of
$$\sup_{y\in{\bf R}^n}\int_{{\bf R}^n}\left|K_{1,m}^\pm(x,y,\lambda_1)-K_{1,m}^\pm(x,y,\lambda_2)\right|^2\langle x\rangle^{-n-2\epsilon}dx$$ $$\le
 C\lambda_1^{n-3}|\lambda_1-\lambda_2|\int_{{\bf R}^n}\langle x\rangle^{-n-2\epsilon}dx\le C\lambda_1^{n-3}|\lambda_1-\lambda_2|.$$
 To prove (3.15) we will use the inequality $$|f(\sigma\lambda_1)-f(\sigma\lambda_2)|\le\sigma\int_{\lambda_1}^{\lambda_2}|f'(\sigma\lambda)|d\lambda\le\sigma|\lambda_1-\lambda_2|^{1/2}\left(\int_{\lambda_1}^{\lambda_2}|f'(\sigma\lambda)|^2d\lambda\right)^{1/2}.\eqno{(3.20)}$$
 Applying (3.20) with $f(z)=\frac{d^m(\phi(z){\cal H}_\nu^\pm(z))}{dz^m}$, $m=\frac{n-2}{2}$, $\sigma=|x-y|$, and using  (3.11), we obtain
$$\left|K_{2,m}^\pm(x,y,\lambda_1)-K_{2,m}^\pm(x,y,\lambda_2)\right|\le C\lambda_1^{-1}|\lambda_1-\lambda_2|^{1/2}|x-y|^{-\frac{n-1}{2}}.\eqno{(3.21)}$$
Hence, the left-hand side of (3.15) is equivalent to
$$\sup_{y\in{\bf R}^n}\int_{{\bf R}^n}\left|K_{2,m}^\pm(x,y,\lambda_1)-K_{2,m}^\pm(x,y,\lambda_2)\right|\langle x\rangle^{-\frac{n+1}{2}-\epsilon}dx$$ $$\le C'\lambda_1^{-1}|\lambda_1-\lambda_2|^{1/2}\sup_{y\in{\bf R}^n}\int_{{\bf R}^n}|x-y|^{-\frac{n-1}{2}}\langle x\rangle^{-\frac{n+1}{2}-\epsilon}dx\le C\lambda_1^{-1}|\lambda_1-\lambda_2|^{1/2}.$$
\eproof

\begin{lemma} There exists a constant $\lambda_0>0$ so that if $n$ is odd, then for all integers $0\le m\le\frac{n-1}{2}$ and all $\lambda\ge \lambda_0$, $0<\epsilon\ll 1$, we have
$$\left\|\langle x\rangle^{-1/2-m-\epsilon}\frac{d^mR_0^\pm(\lambda)}{d\lambda^m}\langle x\rangle^{-1/2-m-\epsilon}\right\|_{L^2\to L^2}$$ $$+
\left\|\langle x\rangle^{-1/2-m-\epsilon}\frac{d^mR^\pm(\lambda)}{d\lambda^m}\langle x\rangle^{-1/2-m-\epsilon}\right\|_{L^2\to L^2}\le C\lambda^{-1}.\eqno{(3.22)}$$
If $n$ is even, (3.22) still holds for $0\le m\le\frac{n-2}{2}$. In this case we also have with $m=\frac{n-2}{2}$
$$\left\|\langle x\rangle^{-1-m-\epsilon}\left(\frac{d^mR_0^\pm}{d\lambda^m}(\lambda_1)-\frac{d^mR_0^\pm}{d\lambda^m}(\lambda_2)\right)\langle x\rangle^{-1-m-\epsilon}\right\|_{L^2\to L^2}$$ $$+
\left\|\langle x\rangle^{-1-m-\epsilon}\left(\frac{d^mR^\pm}{d\lambda^m}(\lambda_1)-\frac{d^mR^\pm}{d\lambda^m}(\lambda_2)\right)\langle x\rangle^{-1-m-\epsilon}\right\|_{L^2\to L^2}\le C\lambda_1^{-1}|\lambda_1-\lambda_2|^{1/2},\eqno{(3.23)}$$
for $\lambda_1+1\ge\lambda_2>\lambda_1\ge \lambda_0$.
\end{lemma}

This lemma is proved in \cite{kn:V} (see Lemma 3.6) and therefore we omit the proof. To prove (3.1) and (3.2) observe first that the operator
$\frac{d^m}{d\lambda^m}T_k^\pm(\lambda)$, $0\le m\le\frac{n-1}{2}$, is a linear combination of operators of the form
$${\cal M}_k^\pm(\lambda,m_1,...,m_{k+1})=\frac{d^{m_1}}{d\lambda^{m_1}}R_0^\pm(\lambda)\left(-V\frac{d^{m_2}}{d\lambda^{m_2}}R_0^\pm(\lambda)\right)...
\left(-V\frac{d^{m_{k+1}}}{d\lambda^{m_{k+1}}}R_0^\pm(\lambda)\right),$$
where $m_j\ge 0$ are integers such that $m_1+...+m_{k+1}\le \frac{n-1}{2}$. Define the operator $\widetilde {\cal M}_k^\pm$ by replacing in the definition of ${\cal M}_k^\pm$ the operator $\frac{d^{m_{k+1}}}{d\lambda^{m_{k+1}}}R_0^\pm(\lambda)$ by $\frac{d^{m_{k+1}}}{d\lambda^{m_{k+1}}}{\cal K}_1^\pm(\lambda)$. We will prove by induction in $k$ that the operators ${\cal M}_k^\pm$ and $\widetilde {\cal M}_k^\pm$ satisfy (3.1) and (3.2). 
Let first $k=1$. Then the kernel, $M_1^\pm$, of the operator ${\cal M}_1^\pm$ satisfies the bound
$$\left|M_1^\pm(x,y,\lambda)\right|$$ $$\le C\int_{{\bf R}^n}|x-\xi|^{m_1-n+2}\left|\left(\partial_z^{m_1}H_\nu^\pm\right)(\lambda|x-\xi|)\right|
|y-\xi|^{m_2-n+2}\left|\left(\partial_z^{m_2}H_\nu^\pm\right)(\lambda|y-\xi|)\right||V(\xi)|d\xi$$
$$\le C\lambda^{n-3}\int_{{\bf R}^n}|x-\xi|^{m_1-\frac{n-1}{2}}
|y-\xi|^{m_2-\frac{n-1}{2}}|V(\xi)|d\xi$$
$$\le C\lambda^{n-3}\int_{{\bf R}^n}\left(|x-\xi|^{m_1+m_2-n+1}+
|y-\xi|^{m_1+m_2-n+1}\right)\langle\xi\rangle^{-\frac{n+1}{2}-\epsilon}d\xi\le C\lambda^{n-3},$$
where we have used (3.10). If $n$ is even, we take $m=\frac{n-2}{2}$ and observe that if $m_1+m_2\le\frac{n-2}{2}$ we have
$$\left|M_1^\pm(x,y,\lambda_1)-M_1^\pm(x,y,\lambda_2)\right|$$ $$\le C\int_{{\bf R}^n}|x-\xi|^{m_1-n+2}\left|\left(\partial_z^{m_1}H_\nu^\pm\right)(\lambda_1|x-\xi|)-\left(\partial_z^{m_1}H_\nu^\pm\right)(\lambda_2|x-\xi|)\right|$$ $$\times
|y-\xi|^{m_2-n+2}\left|\left(\partial_z^{m_2}H_\nu^\pm\right)(\lambda_1|y-\xi|)\right||V(\xi)|d\xi$$
$$+C\int_{{\bf R}^n}|x-\xi|^{m_1-n+2}\left|\left(\partial_z^{m_1}H_\nu^\pm\right)(\lambda_2|x-\xi|)\right|$$ $$\times
|y-\xi|^{m_2-n+2}\left|\left(\partial_z^{m_2}H_\nu^\pm\right)(\lambda_1|y-\xi|)-\left(\partial_z^{m_2}H_\nu^\pm\right)(\lambda_2|y-\xi|)\right||V(\xi)|d\xi$$
$$\le C\lambda_1^{n-3}|\lambda_1-\lambda_2|^{1/2}\int_{{\bf R}^n}\left(|x-\xi|^{m_1-\frac{n-2}{2}}
|y-\xi|^{m_2-\frac{n-1}{2}}+|x-\xi|^{m_1-\frac{n-1}{2}}
|y-\xi|^{m_2-\frac{n-2}{2}}\right)|V(\xi)|d\xi$$
$$\le C\lambda_1^{n-3}|\lambda_1-\lambda_2|^{1/2}\int_{{\bf R}^n}\left(|x-\xi|^{m_1+m_2-n+3/2}+
|y-\xi|^{m_1+m_2-n+3/2}\right)\langle\xi\rangle^{-\frac{n+1}{2}-\epsilon}d\xi$$ $$\le C\lambda_1^{n-3}|\lambda_1-\lambda_2|^{1/2},$$
where we have used (3.20) together with (3.10). Clearly, the operator $\widetilde{\cal M}_1^\pm$
can be treated in precisely the same way. For $k\ge 2$ we have
$$\left\|{\cal M}_k^\pm(\lambda,m_1,...,m_{k+1})\right\|_{L^1\to L^\infty}+\left\|\widetilde{\cal M}_k^\pm(\lambda,m_1,...,m_{k+1})\right\|_{L^1\to L^\infty}$$
$$\le 2\left\|\frac{d^{m_1}}{d\lambda^{m_1}}{\cal K}_1^\pm(\lambda)\left(-V\frac{d^{m_2}}{d\lambda^{m_2}}R_0^\pm(\lambda)\right)...\left(-V\frac{d^{m_k}}{d\lambda^{m_k}}R_0^\pm(\lambda)\right)
\left(-V\frac{d^{m_{k+1}}}{d\lambda^{m_{k+1}}}{\cal K}_1^\pm(\lambda)\right)\right\|_{L^1\to L^\infty}$$
$$+\left\|{\cal M}_{k-1}^\pm(\lambda,m_1,...,m_{k})\right\|_{L^1\to L^\infty}\left\|V\frac{d^{m_{k+1}}}{d\lambda^{m_{k+1}}}{\cal K}_2^\pm(\lambda)\right\|_{L^1\to L^1}$$
$$+2\left\|\frac{d^{m_{1}}}{d\lambda^{m_{1}}}{\cal K}_2^\pm(\lambda)V\right\|_{L^\infty\to L^\infty}\left\|\widetilde{\cal M}_{k-1}^\pm(\lambda,m_2,...,m_{k+1})\right\|_{L^1\to L^\infty}$$
 $$\le C\left\|\frac{d^{m_1}}{d\lambda^{m_1}}{\cal K}_1^\pm(\lambda)\langle x\rangle^{-1/2-m_1-\epsilon}\right\|_{L^2\to L^\infty}\left\|\langle x\rangle^{-1/2-m_2-\epsilon}\frac{d^{m_2}}{d\lambda^{m_2}}R_0^\pm(\lambda)\langle x\rangle^{-1/2-m_2-\epsilon}\right\|_{L^2\to L^2}...$$ $$\times\left\|\langle x\rangle^{-1/2-m_k-\epsilon}\frac{d^{m_k}}{d\lambda^{m_k}}R_0^\pm(\lambda)\langle x\rangle^{-1/2-m_k-\epsilon}\right\|_{L^2\to L^2}\left\|\langle x\rangle^{-1/2-m_{k+1}-\epsilon}
\frac{d^{m_{k+1}}}{d\lambda^{m_{k+1}}}{\cal K}_1^\pm(\lambda)\right\|_{L^1\to L^2}$$
$$+C\left\|{\cal M}_{k-1}^\pm(\lambda,m_1,...,m_{k})\right\|_{L^1\to L^\infty}\left\|\langle x\rangle^{-1-m_{k+1}-\epsilon}\frac{d^{m_{k+1}}}{d\lambda^{m_{k+1}}}{\cal K}_2^\pm(\lambda)\right\|_{L^1\to L^1}$$
$$+C\left\|\frac{d^{m_{1}}}{d\lambda^{m_{1}}}{\cal K}_2^\pm(\lambda)\langle x\rangle^{-1-m_{1}-\epsilon}\right\|_{L^\infty\to L^\infty}\left\|\widetilde{\cal M}_{k-1}^\pm(\lambda,m_2,...,m_{k+1})\right\|_{L^1\to L^\infty}$$
$$\le C\lambda^{n-2-k}+C\lambda^{-1}\left\|{\cal M}_{k-1}^\pm(\lambda,m_1,...,m_{k})\right\|_{L^1\to L^\infty}
+C\lambda^{-1}\left\|\widetilde{\cal M}_{k-1}^\pm(\lambda,m_2,...,m_{k+1})\right\|_{L^1\to L^\infty},$$
where we have used Lemmas 3.2 and 3.3. Therefore, if (3.1) is valid for $k-1$ it is also valid for $k$.
The bound (3.2) can be proved in precisely the same way using (3.14), (3.15) and (3.23). Furthermore, using the resolvent identity
$$R^\pm(\lambda)-R_0^\pm(\lambda)=-R_0^\pm(\lambda)VR^\pm(\lambda)=-R^\pm(\lambda)VR_0^\pm(\lambda)$$
together with (3.1) and proceeding as above, one easily obtains
$$\sum_{m_1+...m_{k+1}\le\frac{n-1}{2}}\left\|\frac{d^{m_1}}{d\lambda^{m_1}}R^\pm(\lambda)\left(-V\frac{d^{m_2}}{d\lambda^{m_2}}R_0^\pm(\lambda)\right)...
\left(-V\frac{d^{m_{k+1}}}{d\lambda^{m_{k+1}}}R_0^\pm(\lambda)\right)\right\|_{L^1\to L^\infty}\le C\lambda^{n-2-k}$$ $$+C\lambda^{-1}\sum_{m_1+...m_{k+1}\le\frac{n-1}{2}}\left\|\frac{d^{m_1}}{d\lambda^{m_1}}R^\pm(\lambda)\left(-V\frac{d^{m_2}}{d\lambda^{m_2}}R_0^\pm(\lambda)\right)...
\left(-V\frac{d^{m_{k+1}}}{d\lambda^{m_{k+1}}}R_0^\pm(\lambda)\right)\right\|_{L^1\to L^\infty}.\eqno{(3.24)}$$
Taking $\lambda$ big enough one can absorbe the second term in the right-hand side of (3.24) and conclude that the operator $\widetilde T_k^\pm$ 
satisfies (3.1), too. Similarly, it is easy to see that $\widetilde T_k^\pm$ satisfies also (3.2).
\eproof

\section{Study of the operator ${\cal A}_1$}

In this section we will prove the following

\begin{Theorem} Suppose that $V$ satisfies (1.5). Then in all dimensions $n\ge 4$ we have the estimate
$$\left\|{\cal A}_1^\pm(t,h)\right\|_{L^1\to L^\infty}\le Ch^{-\frac{n+1}{2}}|t|^{-\frac{n-1}{2}}.\eqno{(4.1)}$$
\end{Theorem}

{\it Proof.} Clearly, it suffices to consider the case $"-"$  and $t>0$, only. It is easy to see that 
the kernel of the operator ${\cal A}_1^-$ is of the form
$$\int_{{\bf R}^n}A_h(|x-\xi|,|y-\xi|,t)V(\xi)d\xi,$$
where
$$A_h(\sigma_1,\sigma_2,t)=\frac{(\sigma_1\sigma_2)^{2-n}}{ih2^3(2\pi)^{n-1}}\int e^{it\lambda}\widetilde\varphi(h\lambda){\cal H}_\nu^-(\sigma_1\lambda){\cal H}_\nu^-(\sigma_2\lambda)d\lambda.$$
Let $\phi\in C_0^\infty({\bf R})$, $\phi(z)=1$ for $|z|\le 1/2$, $\phi(z)=0$ for $|z|\ge 1$. Decompose the function $A_h$ as
$A_h^{(1)}+A_h^{(2)}$, where
$$A_h^{(1)}(\sigma_1,\sigma_2,t)=\frac{(\sigma_1\sigma_2)^{2-n}}{ih2^3(2\pi)^{n-1}}\int e^{it\lambda}\widetilde\varphi(h\lambda)((1-\phi){\cal H}_\nu^-)(\sigma_1\lambda)((1-\phi){\cal H}_\nu^-)(\sigma_2\lambda)d\lambda,$$
$$A_h^{(2)}(\sigma_1,\sigma_2,t)=$$ $$\frac{(\sigma_1\sigma_2)^{2-n}}{ih2^3(2\pi)^{n-1}}\int e^{it\lambda}\widetilde\varphi(h\lambda)\left((\phi {\cal H}_\nu^-)(\sigma_1\lambda){\cal H}_\nu^-(\sigma_2\lambda)+((1-\phi){\cal H}_\nu^-)(\sigma_1\lambda)(\phi {\cal H}_\nu^-)(\sigma_2\lambda)\right)d\lambda.$$
Thus we decompose the operator ${\cal A}_1^-$ as ${\cal A}^{(1)}+{\cal A}^{(2)}$, where ${\cal A}^{(j)}$ is defined by replacing in the definition of ${\cal A}_1^-$ the function $A_h$ by $A_h^{(j)}$. 

\begin{lemma} The operator ${\cal A}^{(2)}$ satisfies the estimate
$$\left\|{\cal A}^{(2)}(t,h)\right\|_{L^1\to L^\infty}\le Ch^{-\frac{n-1}{2}}t^{-\frac{n-1}{2}}.\eqno{(4.2)}$$
\end{lemma}

{\it Proof.} Clearly, it suffices to show that the function $A_h^{(2)}$ satisfies the bound
$$\left|A_h^{(2)}(\sigma_1,\sigma_2,t)\right|\le Ch^{-\frac{n-1}{2}}t^{-\frac{n-1}{2}}\left(\sigma_1^{-n+2}+\sigma_2^{-n+2}\right).\eqno{(4.3)}$$
Since
$$A_h^{(2)}(\sigma_1,\sigma_2,t)=h^{-2n+2}A_1^{(2)}\left(\frac{\sigma_1}{h},\frac{\sigma_2}{h},\frac{t}{h}\right),$$
it suffices to prove (4.3) for $h=1$. To do so, recall that the function ${\cal H}_\nu^-$ satisfies (3.10) and (3.11). 
Hence, for $\lambda\in{\rm supp}\,\widetilde\varphi$ and all integers $0\le k\le\frac{n}{2}$, we have
$$\left|\frac{d^k}{d\lambda^k}{\cal H}_\nu^-(\sigma\lambda)\right|\le C\sigma^{\frac{n-3}{2}+k},\quad\forall \sigma>0,\eqno{(4.4)}$$
 $$\left|\frac{d^k}{d\lambda^k}\left(\phi{\cal H}_\nu^-\right)(\sigma\lambda)\right|\le C,\quad\forall \sigma>0.\eqno{(4.5)}$$
Let $0\le m\le\frac{n}{2}$ be an integer. 
Integrating by parts $m$ times and using (4.4) and (4.5) we conclude that the function $A_1^{(2)}$ satisfies
the bound
 $$\left|A_1^{(2)}(\sigma_1,\sigma_2,t)\right|\le Ct^{-m}(\sigma_1\sigma_2)^{-n+2}\left(\sigma_1^{m+\frac{n-3}{2}}+
 \sigma_2^{m+\frac{n-3}{2}}\right),\eqno{(4.6)}$$
for all integers $0\le m\le\left[\frac{n}{2}\right]$ and hence, by interpolation, for all real
$0\le m\le\left[\frac{n}{2}\right]$. Taking $m=\frac{n-1}{2}$ in (4.6) we get (4.3) with $h=1$.
\eproof

\begin{prop} The operator ${\cal A}^{(1)}$ satisfies the estimate
$$\left\|{\cal A}^{(1)}(t,h)\right\|_{L^1\to L^\infty}\le Ch^{-\frac{n+1}{2}}t^{-\frac{n-1}{2}}.\eqno{(4.7)}$$
\end{prop}

{\it Proof.}  Recall first that ${\cal H}_\nu^-(z)=e^{-iz}b_\nu^-(z)$, where $b_\nu^-(z)$ is a symbol of order $\frac{n-3}{2}$. We would like to integrate by parts $\frac{n-3}{2}$ times with respect to the variable $\xi$. Set
$$\psi(x,\xi,y)=|x-\xi|+|y-\xi|,\quad \rho(x,\xi,y)=\nabla_\xi \psi(x,\xi,y)=\frac{\xi -x}{|\xi -x|}+\frac{\xi -y}{|\xi -y|}.$$
We are going to use the identity
$$e^{-i\lambda\psi}=\Lambda_\xi e^{-i\lambda\psi},$$
where
$$\Lambda_\xi=\frac{i\rho}{\lambda|\rho|^2}\cdot\nabla_\xi.$$
Clearly, the function $\rho$ satisfies the bound
$$\left|\partial_\xi^\alpha\rho(x,\xi,y)\right|\le C_\alpha\left(|x-\xi|^{-|\alpha|}+|y-\xi|^{-|\alpha|}\right).\eqno{(4.8)}$$
Using (4.8) one can easily see by induction that the operator $(\Lambda_\xi^*)^m$ is of the form
$$\lambda^{-m}\sum_{0\le|\alpha|\le m} r_\alpha^{(m)}(x,\xi,y)\partial_\xi^\alpha$$
with functions $r_\alpha^{(m)}$ satisfying the bound
$$\left|r_\alpha^{(m)}(x,\xi,y)\right|\le C_m|\rho|^{-2m+|\alpha|}\left(|x-\xi|^{-m+|\alpha|}+|y-\xi|^{-m+|\alpha|}\right).\eqno{(4.9)}$$
Observe also that making a change of variables $\xi\to h\xi$ we can write the kernel of the operator ${\cal A}^{(1)}(t,h)$ in the form
$$h^{-n+2}\int_{{\bf R}^n}A^{(1)}_1(|x'-\xi|,|y'-\xi|,t')V(h\xi)d\xi,\eqno{(4.10)}$$
where $x'=x/h$, $y'=y/h$, $t'=t/h$. Let $1\le m\le\frac{n-1}{2}$ be an integer. Integrating by parts $m$ times with respect to the variable $\lambda$ we can write the function $A_1^{(1)}$ in the form
$$\sum_{j_1+j_2\le m}A_{1,j_1,j_2}^{(1)},$$
where 
$$A_{1,j_1,j_2}^{(1)}(\sigma_1,\sigma_2,t)=t^{-m}\int e^{i(t-(\sigma_1+\sigma_2))\lambda}\varphi_{j_1,j_2}(\lambda)b_{\nu,j_1}^-(\sigma_1\lambda)b_{\nu,j_2}^-(\sigma_2\lambda)d\lambda,$$
 $\varphi_{j_1,j_2}\in C_0^\infty((0,+\infty))$ and 
 $$b_{\nu,j}^-(z)=z^{j-n+2}e^{iz}\frac{d^j}{dz^j}\left(e^{-iz}(1-\phi(z))b_\nu^-(z)\right).$$
 Since $b_\nu^-(z)$ is a symbol of order $\frac{n-3}{2}$ for $z\ge 1$, we have the bound
 $$\left|\frac{d^k}{dz^k}b_{\nu,j}^-(z)\right|\le C_{j,k}\langle z\rangle^{j-\frac{n-1}{2}-k},\quad\forall z>0.\eqno{(4.11)}$$
 Denote by ${\cal A}_{j_1,j_2}^{(1)}(t,h)$ the operator with kernel
$$h^{-n+2}\int_{{\bf R}^n}A^{(1)}_{1,j_1,j_2}(|x'-\xi|,|y'-\xi|,t')V(h\xi)d\xi.\eqno{(4.12)}$$
Integrating by parts $m-1$ times with respect to $\xi$ in the integral in (4.12)
we write the kernel of the operator ${\cal A}_{j_1,j_2}^{(1)}(t,h)$ in the form
$$(t/h)^{-m}h^{-n+2}\sum_{0\le|\alpha|\le m-1}\int\int_{{\bf R}^n}e^{i\lambda(t'-|x'-\xi|-|y'-\xi|)}
 \lambda^{-m+1}\varphi_{j_1,j_2}(\lambda)r_\alpha^{(m-1)}(x',\xi,y')$$ 
$$\times\sum_{|\alpha_1|+|\alpha_2|=|\alpha|}c_{\alpha,\alpha_1,\alpha_2}R_{\alpha_1}(x',\xi,y',\lambda)\partial_\xi^{\alpha_2}\left(V(h\xi)\right)d\xi d\lambda,$$
 where $c_{\alpha,\alpha_1,\alpha_2}$ are constants and  
$$R_\alpha=\partial_\xi^\alpha\left(b_{\nu,j_1}^-(\lambda|x'-\xi|)b_{\nu,j_2}^-(\lambda|y'-\xi|)\right).$$
It follows from (4.11) that as long as $\lambda\in{\rm supp}\,\varphi_{j_1,j_2}$ we have the bound
 $$\left|R_{\alpha}(x',\xi,y',\lambda)\right|\le C_\alpha\langle x'-\xi\rangle^{m+1-n-|\alpha|}+C_\alpha\langle y'-\xi\rangle^{m+1-n-|\alpha|}.\eqno{(4.13)}$$
 By (4.9) and (4.13) we obtain
 $$(t/h)^mh^{n-2}\left\|{\cal A}_{j_1,j_2}^{(1)}(t,h)\right\|_{L^1\to L^\infty}\le 
  $$ $$C\sum_{0\le|\alpha|\le m-1}\int_{{\bf R}^n}|\rho(x',\xi,y')|^{-2m+2+|\alpha|}\left(| x'-\xi|^{-m+1+|\alpha|}+|y'-\xi|^{-m+1+|\alpha|}\right)$$ 
  $$\times \sum_{|\alpha_1|+|\alpha_2|=|\alpha|}\left(|x'-\xi|^{m+1-n-|\alpha_1|}+| y'-\xi|^{m+1-n-|\alpha_1|}\right)\left|\partial_\xi^{\alpha_2}\left(V(h\xi)\right)\right|d\xi $$ $$\le 
  C\sum_{0\le|\alpha_2|\le m-1}\int_{{\bf R}^n}|\rho(x',\xi,y')|^{-2m+2+|\alpha_2|}\left(| x'-\xi|^{-n+2+|\alpha_2|}+| y'-\xi|^{-n+2+|\alpha_2|}\right)\left|\partial_\xi^{\alpha_2}\left(V(h\xi)\right)\right|d\xi $$
  $$\le Ch^{-2}
  \sum_{0\le|\alpha_2|\le m-1}\int_{{\bf R}^n}|\rho(x,\xi,y)|^{-2m+2+|\alpha_2|}\left(|x-\xi|^{-n+2+|\alpha_2|}+| y-\xi|^{-n+2+|\alpha_2|}\right)\left|\partial_\xi^{\alpha_2}V(\xi)\right|d\xi $$
  $$\le Ch^{-2}\|V\|_{{\cal V}_{2+\epsilon'}^{m-1}}$$ $$\times 
  \sum_{0\le|\alpha_2|\le m-1}\int_{{\bf R}^n}|\rho(x,\xi,y)|^{-2m+2+|\alpha_2|}\left(|x-\xi|^{-n+2+|\alpha_2|}+| y-\xi|^{-n+2+|\alpha_2|}\right)\langle\xi\rangle^{-2-\epsilon'-|\alpha_2|}d\xi. $$
  We need now the following bound proved in the appendix.
  
  \begin{lemma} Let $0\le\ell_1< n-1$, $\ell_1\le\ell_2< n$, $\ell_2+\ell_3>n$. Then we have the bound
  $$\sup_{x,y\in{\bf R}^n}\int_{{\bf R}^n}|\rho(x,\xi,y)|^{-\ell_1}|x-\xi|^{-\ell_2}\langle\xi\rangle^{-\ell_3}d\xi<+\infty.\eqno{(4.14)}$$
  \end{lemma}
 
  Thus we conclude that
  $$\left\|{\cal A}^{(1)}(t,h)\right\|_{L^1\to L^\infty}\le C\|V\|_{{\cal V}_{2+\epsilon'}^{m-1}}h^{m-n}t^{-m}.\eqno{(4.15)}$$
  If $n$ is odd (4.7) follows from (4.15) applied with $m=\frac{n-1}{2}$. Let now $n$ be even. Then there exists a family of functions, $V_\theta$, $0<\theta\le 1$, such that
  $$\left\|V_\theta\right\|_{{\cal V}_{2+\epsilon'}^{\frac{n-4}{2}}}+\theta^{-1/2}\left\|V-V_\theta\right\|_{{\cal V}_{2+\epsilon'}^{\frac{n-4}{2}}}+\theta^{1/2}\left\|V_\theta\right\|_{{\cal V}_{2+\epsilon'}^{\frac{n-2}{2}}}\le C.\eqno{(4.16)}$$
  To prove the desired estimate in this case we are going to apply interpolation arguments. To this end, we will make use of the partition of the unity
  $$1=\sum_{k=0}^\infty\phi_k(\lambda),\quad \lambda>0,\eqno{(4.17)}$$
  where $\phi_0=\phi$ (the function $\phi$ being defined above), $\phi_k(\lambda)=\phi^\sharp(2^{-k}\lambda)$, $k\ge 1$, with a function $\phi^\sharp\in C_0^\infty({\bf R})$, $\phi^\sharp(\lambda)=0$ for $0<\lambda\le 1/2$ and $\lambda\ge 1$. It follows from (4.16) that the functions
  $$V_k^{(0)}(\xi)=(V(\xi)-V_\theta(\xi))\phi_k(\langle\xi\rangle),\quad V_k^{(1)}(\xi)=V_\theta(\xi)\phi_k(\langle\xi\rangle),$$
  satisfy the bound
   $$\left\|V_k^{(j)}\right\|_{{\cal V}_{5/2-j+\epsilon'/2}^{\frac{n-4}{2}+j}}\le C2^{-k\epsilon'/2}\left(2^k\theta\right)^{1/2-j}.\eqno{(4.18)}$$
  We now decompose the function (4.10) as
   $$ \sum_{k=0}^\infty \sum_{j=0}^1 F_k^{(j)},$$ where 
   $F_k^{(j)}$ is defined by replacing in (4.10) the function $V$ by $V_k^{(j)}$. 
Clearly, it suffices to show that 
  $$\left|\sum_{j=0}^1F_k^{(j)}\right|\le C2^{-k\epsilon_0} h^{-\frac{n+1}{2}}t^{-\frac{n-1}{2}},\eqno{(4.19)}$$
  with some constant $\epsilon_0>0$. This in turn follows from the following estimates applied with $\theta=h2^{-k}$.
  
  \begin{lemma} The functions $F_k^{(j)}$, $j=0,1$, satisfy the estimate
  $$\left|F_k^{(j)}\right|\le C2^{-k\epsilon_0} h^{-\frac{n+1}{2}}t^{-\frac{n-1}{2}}\left(\frac{\theta 2^{k}}{h}\right)^{1/2-j},\eqno{(4.20)}$$
   with some $\epsilon_0>0$.
  \end{lemma}
  
  {\it Proof.} We first integrate by parts $\frac{n-4}{2}+j$ times with respect to the variable $\xi$. Thus we get that $F_k^{(j)}$ is a linear combination of functions of the form
  $$h^{-n+2}\int_{{\bf R}^n}f_{\alpha_1,\alpha_2}^{(j)}(t',x',y',\xi)r_\alpha^{(\frac{n-4}{2}+j)}(x',\xi,y')\partial_\xi^\beta\left( V_k^{(j)}(h\xi)\right)d\xi,$$
  where
  $$f_{\alpha_1,\alpha_2}^{(j)}=\int e^{i\lambda(t'-|x'-\xi|-|y'-\xi|)}\lambda^{-\frac{n-4}{2}-j}
  \widetilde\varphi(\lambda)\partial_\xi^{\alpha_1}\left(b^-_{\nu,0}(\lambda|x'-\xi|)\right)
  \partial_\xi^{\alpha_2}\left(b^-_{\nu,0}(\lambda|y'-\xi|)\right)d\lambda,$$
  $|\alpha_1|+|\alpha_2|+|\beta|=|\alpha|\le \frac{n-4}{2}+j$. Let $0\le m\le \frac{n}{2}$ be an integer. 
  We now integrate by parts $m$ times with respect to $\lambda$ to obtain
  $$\left|f_{\alpha_1,\alpha_2}^{(j)}\right|\le C(t')^{-m}\sum_{j_1+j_2\le m}\int_{{\rm supp}\,\widetilde\varphi}\left|\partial_\xi^{\alpha_1}\left(b^-_{\nu,j_1}(\lambda|x'-\xi|)\right)\right|\left|\partial_\xi^{\alpha_2}\left(b^-_{\nu,j_2}(\lambda|y'-\xi|)\right)\right|d\lambda$$
   $$\le C(t')^{-m}\sum_{j_1+j_2\le m}\langle x'-\xi\rangle^{j_1-\frac{n-1}{2}-|\alpha_1|}\langle y'-\xi\rangle^{j_2-\frac{n-1}{2}-|\alpha_2|}$$
  $$\le C(t')^{-m}\left(\langle x'-\xi\rangle+\langle y'-\xi\rangle\right)^m\langle x'-\xi\rangle^{-\frac{n-1}{2}-|\alpha_1|}\langle y'-\xi\rangle^{-\frac{n-1}{2}-|\alpha_2|},
  \eqno{(4.21)}$$
  where we have used (4.11). By interpolation, (4.21) holds for all real $0\le m\le \frac{n}{2}$ and in particular for $m=\frac{n-1}{2}$. Hence, using this together with (4.9), we obtain
  $$\left|f_{\alpha_1,\alpha_2}^{(j)}r_\alpha^{(\frac{n-4}{2}+j)}\right|\le C(t')^{-\frac{n-1}{2}}\left(|x'-\xi|^{-1}+|y'-\xi|^{-1}\right)^{\frac{n-1}{2}+|\alpha_1|+|\alpha_2|}\left|r_\alpha^{(\frac{n-4}{2}+j)}(x',\xi,y')
  \right|$$ $$\le C(t')^{-\frac{n-1}{2}}|\rho|^{-n+4-2j+|\beta|}\left(|x'-\xi|^{-1}+| y'-\xi|^{-1}\right)^{n-\frac{5}{2}+j-|\beta|}.\eqno{(4.22)}$$
  By (4.18) and (4.22), we obtain
  $$(t/h)^{\frac{n-1}{2}}h^{n}\left|F_k^{(j)}\right|\le C2^{-k\epsilon'/2}\left(2^k\theta/h\right)^{1/2-j}
  \sum_{0\le|\beta|\le \frac{n-4}{2}+j}\int_{{\bf R}^n}|\rho(x,\xi,y)|^{-n+4-2j+|\beta|}$$
  $$\times \left(|x-\xi|^{-n+5/2-j+|\beta|}+| y-\xi|^{-n+5/2-j+|\beta|}\right)\langle\xi\rangle^{-5/2+j-\epsilon'/2-|\beta|}d\xi,$$
  which together with Lemma 4.4 imply (4.20).
  \eproof

\section{Study of the operator ${\cal A}_2$}

In this section we will prove the following

\begin{Theorem} Suppose that $V$ satisfies (1.6). Then in dimensions $n=6,7$ we have the estimate 
$$\left\|{\cal A}_2^\pm(t,h)\right\|_{L^1\to L^\infty}\le C_\epsilon h^{-\frac{n+1}{2}-\epsilon}|t|^{-\frac{n-1}{2}},\eqno{(5.1)}$$
for every $0<\epsilon\ll 1$.
\end{Theorem}

{\it Proof.} Clearly, it suffices to consider the case $"-"$  and $t>0$, only. The kernel of the operator ${\cal A}_2^-$ is of the form
$$\int_{{\bf R}^n}\int_{{\bf R}^n}B_h(|x-\xi_1|,|\xi_1-\xi_2|,|\xi_2-y|,t)V(\xi_1)V(\xi_2)d\xi_1d\xi_2,$$
where
$$B_h(\sigma_1,\sigma_2,\sigma_3,t)=\frac{(\sigma_1\sigma_2\sigma_3)^{2-n}}{h2^5(2\pi)^{\frac{3n}{2}-2}}\int e^{it\lambda}\widetilde\varphi(h\lambda){\cal H}_\nu^-(\sigma_1\lambda){\cal H}_\nu^-(\sigma_2\lambda){\cal H}_\nu^-(\sigma_3\lambda)d\lambda.$$
As in the previous section, we decompose the function $B_h$ as $\sum_{j=1}^4B_h^{(j)}$, where
$$B_h^{(1)}(\sigma_1,\sigma_2,\sigma_3,t)=$$ $$\frac{(\sigma_1\sigma_2\sigma_3)^{2-n}}{h2^5(2\pi)^{\frac{3n}{2}-2}}\int e^{it\lambda}\widetilde\varphi(h\lambda)((1-\phi){\cal H}_\nu^-)(\sigma_1\lambda)((1-\phi){\cal H}_\nu^-)(\sigma_2\lambda)((1-\phi){\cal H}_\nu^-)(\sigma_3\lambda)d\lambda,$$
$$B_h^{(2)}(\sigma_1,\sigma_2,\sigma_3,t)=$$ $$\frac{(\sigma_1\sigma_2\sigma_3)^{2-n}}{h2^5(2\pi)^{\frac{3n}{2}-2}}\int e^{it\lambda}\widetilde\varphi(h\lambda)(\phi{\cal H}_\nu^-)(\sigma_1\lambda){\cal H}_\nu^-(\sigma_2\lambda){\cal H}_\nu^-(\sigma_3\lambda)d\lambda,$$
$$B_h^{(3)}(\sigma_1,\sigma_2,\sigma_3,t)=$$ $$\frac{(\sigma_1\sigma_2\sigma_3)^{2-n}}{h2^5(2\pi)^{\frac{3n}{2}-2}}\int e^{it\lambda}\widetilde\varphi(h\lambda)((1-\phi){\cal H}_\nu^-)(\sigma_1\lambda)(\phi{\cal H}_\nu^-)(\sigma_2\lambda){\cal H}_\nu^-(\sigma_3\lambda)d\lambda,$$
$$B_h^{(4)}(\sigma_1,\sigma_2,\sigma_3,t)=$$ $$\frac{(\sigma_1\sigma_2\sigma_3)^{2-n}}{h2^5(2\pi)^{\frac{3n}{2}-2}}\int e^{it\lambda}\widetilde\varphi(h\lambda)((1-\phi){\cal H}_\nu^-)(\sigma_1\lambda)((1-\phi){\cal H}_\nu^-)(\sigma_2\lambda)(\phi{\cal H}_\nu^-)(\sigma_3\lambda)d\lambda.$$
Thus we decompose the operator ${\cal A}_2^-$ as $\sum_{j=1}^4{\cal B}^{(j)}$, where ${\cal B}^{(j)}$ is defined by replacing in the definition of ${\cal A}_2^-$ the function $B_h$ by $B_h^{(j)}$. 

\begin{prop} Let $V$ satisfy (1.6). Then in all dimensions $n\ge 6$ the operators ${\cal B}^{(j)}$, $j=2,3,4$, satisfy the estimate
$$\left\|{\cal B}^{(j)}(t,h)\right\|_{L^1\to L^\infty}\le Ch^{-\frac{n+1}{2}}t^{-\frac{n-1}{2}}.\eqno{(5.2)}$$
\end{prop}

{\it Proof.} The operators ${\cal B}^{(j)}$, $j=2,4$, can be treated in precisely the same way as the operator ${\cal A}_1^-$ in the previous section, integrating by parts with respect to the variables $\xi_2$ and $\xi_1$, respectively. Decompose the function $B_h^{(3)}$ as $ B_h^{(3,1)}+ B_h^{(3,2)}$, where $ B_h^{(3,j)}$, $j=1,2$,
are defined by replacing in the definition of $B_h^{(3)}$ the function ${\cal H}_\nu^-(\sigma_3\lambda)$ by $((1-\phi){\cal H}_\nu^-)(\sigma_3\lambda)$ and $(\phi{\cal H}_\nu^-)(\sigma_3\lambda)$, repsectively. Define the operators ${\cal B}_j^{(3)}$, $j=1,2$, by replacing in the definition of ${\cal B}^{(3)}$ the function $B_h^{(3)}$ by $B_h^{(3,j)}$. 
The operator ${\cal B}^{(3)}_2$ can be treated in the same way as the operator ${\cal A}^{(2)}$ in the previous section with no need of integrating by parts with respect to the variables $\xi_1$ and $\xi_2$. 
The analysis of the operator ${\cal B}^{(3)}_1$, 
however, is more complicated and cannot be carried out as in the case of the other operators above. Indeed, to avoid non-integrable singularities at $\xi_1-\xi_2=0$ one needs to proceed differently. The idea is inspired from \cite{kn:EG} and consists of integrating by parts with respect to the variable $\xi_1+\xi_2$ and using the fact that given any smooth function $f$ on ${\bf R}^n$ we have the identity
$$\frac{\partial f(\xi_1-\xi_2)}{\partial(\xi_1+\xi_2)}=0.$$
Set
$$\Psi(x,\xi_1,\xi_2,y)=|x-\xi_1|+|y-\xi_2|,$$ $$ \mu(x,\xi_1,\xi_2,y)=\nabla_{\xi_1+\xi_2} \Psi(x,\xi_1,\xi_2,y)=\frac{\xi_1 -x}{|\xi_1 -x|}+\frac{\xi_2 -y}{|\xi_2 -y|}.$$
We are going to use the identity
$$e^{-i\lambda\Psi}={\cal L}_{\xi_1+\xi_2} e^{-i\lambda\Psi},$$
where
$${\cal L}_{\xi_1+\xi_2}=\frac{i\mu}{\lambda|\mu|^2}\cdot\nabla_{\xi_1+\xi_2}.$$
Clearly, the function $\mu$ satisfies the bound
$$\left|\partial_{\xi_1+\xi_2}^\alpha\mu(x,\xi_1,\xi_2,y)\right|\le C_\alpha\left(|x-\xi_1|^{-|\alpha|}+|y-\xi_2|^{-|\alpha|}\right).\eqno{(5.3)}$$
Using (5.3) one can easily see by induction that the operator $({\cal L}_{\xi_1+\xi_2}^*)^m$ is of the form
$$\lambda^{-m}\sum_{0\le|\alpha|\le m} q_\alpha^{(m)}(x,\xi_1,\xi_2,y)\partial_{\xi_1+\xi_2}^\alpha$$
with functions $q_\alpha^{(m)}$ satisfying the bound
$$\left|q_\alpha^{(m)}(x,\xi_1,\xi_2,y)\right|\le C_m|\mu|^{-2m+|\alpha|}\left(|x-\xi_1|^{-m+|\alpha|}+|y-\xi_2|^{-m+|\alpha|}\right).\eqno{(5.4)}$$
Making a change of variables $\xi_1\to h\xi_1$, $\xi_2\to h\xi_2$, we can write the kernel of the operator ${\cal B}^{(3)}_1(t,h)$ in the form
$$h^{-n+4}\int_{{\bf R}^n}\int_{{\bf R}^n}B_1^{(3,1)}(|x'-\xi_1|,|\xi_1-\xi_2|,|\xi_2-y'|,t')V(h\xi_1)V(h\xi_2)d\xi_1d\xi_2,\eqno{(5.5)}$$
where $x'=x/h$, $y'=y/h$, $t'=t/h$.  Let $1\le m\le\frac{n-1}{2}$ be an integer. Integrating by parts $m$ times with respect to the variable $\lambda$ we can write the function $B_1^{(3,1)}$ in the form
$$\sum_{j_1+j_2+j_3\le m}B^{(3,1)}_{1,j_1,j_2,j_3},$$
where 
$$B^{(3,1)}_{1,j_1,j_2,j_3}(\sigma_1,\sigma_2,\sigma_3,t)=t^{-m}\int e^{i(t-(\sigma_1+\sigma_3))\lambda}\varphi_{j_1,j_2,j_3}(\lambda)b_{\nu,j_1}^-(\sigma_1\lambda)a_{\nu,j_2}^-(\sigma_2\lambda)b_{\nu,j_3}^-(\sigma_3\lambda)d\lambda,$$
 $\varphi_{j_1,j_2,j_3}\in C_0^\infty((0,+\infty))$, $b_{\nu,j}^-(z)$ is as in the previous section, and 
$$a_{\nu,j}^-(z)=z^{j-n+2}\frac{d^j}{dz^j}(\phi{\cal H}_\nu^-)(z).$$
It is easy to see that the function $a_{\nu,j}^-$ satisfies the bound
 $$\left|a_{\nu,j}^-(z)\right|\le Cz^{-n+2},\quad\forall z>0.\eqno{(5.6)}$$
 Denote by ${\cal B}_{j_1,j_2,j_3}^{(3,1)}(t,h)$ the operator with kernel
 $$h^{-n+4}\int_{{\bf R}^n}\int_{{\bf R}^n}B_{1,j_1,j_2,j_3}^{(3,1)}(|x'-\xi_1|,|\xi_1-\xi_2|,|\xi_2-y'|,t')V(h\xi_1)V(h\xi_2)d\xi_1d\xi_2.\eqno{(5.7)}$$
Integrating by parts $m-1$ times with respect to $\xi_1+\xi_2$ in the integral in (5.7)
we obtain
$$(t/h)^{-m}h^{-n+4}\sum_{0\le|\alpha|\le m-1}\int\int_{{\bf R}^n}\int_{{\bf R}^n}e^{i\lambda(t'-|x'-\xi_1|-|y'-\xi_2|)}
 \lambda^{-m+1}\varphi_{j_1,j_2,j_3}(\lambda)a^-_{\nu,j_2}(\lambda|\xi_1-\xi_2|)$$ 
$$\times q_\alpha^{(m-1)}(x',\xi_1,\xi_2,y')\sum_{|\alpha_1|+|\alpha_2|=|\alpha|}c_{\alpha,\alpha_1,\alpha_2}Q_{\alpha_1}(x',\xi_1,\xi_2,y',\lambda)\partial_{\xi_1+\xi_2}^{\alpha_2}\left(V(h\xi_1)V(h\xi_2)\right)d\xi_1d\xi_2 d\lambda,$$
where $c_{\alpha,\alpha_1,\alpha_2}$ are constants and  
$$Q_\alpha=\partial_{\xi_1+\xi_2}^\alpha\left(b_{\nu,j_1}^-(\lambda|x'-\xi_1|)b_{\nu,j_3}^-(\lambda|y'-\xi_2|)\right).$$
It follows from (4.11) that as long as $\lambda\in{\rm supp}\,\varphi_{j_1,j_2,j_3}$ we have the bound
 $$\left|Q_{\alpha}(x',\xi_1,\xi_2,y',\lambda)\right|\le C_\alpha\langle x'-\xi_1\rangle^{m+1-n-|\alpha|}+C_\alpha\langle y'-\xi_2\rangle^{m+1-n-|\alpha|}.\eqno{(5.8)}$$
 By (5.6) and (5.8) we obtain
 $$(t/h)^mh^{n-4}\left\|{\cal B}_{j_1,j_2,j_3}^{(3,1)}(t,h)\right\|_{L^1\to L^\infty}\le 
  $$ $$C\sum_{0\le|\alpha|\le m-1}\int_{{\bf R}^n}\int_{{\bf R}^n}|\mu(x',\xi_1,\xi_2,y')|^{-2m+2+|\alpha|}\left(| x'-\xi_1|^{-m+1+|\alpha|}+|y'-\xi_2|^{-m+1+|\alpha|}\right)$$ 
  $$\times \sum_{|\alpha_1|+|\alpha_2|=|\alpha|}\left(|x'-\xi_1|^{m+1-n-|\alpha_1|}+| y'-\xi_2|^{m+1-n-|\alpha_1|}\right)\left|\partial_{\xi_1+\xi_2}^{\alpha_2}\left(V(h\xi_1)V(h\xi_2)\right)\right|\frac{d\xi_1d\xi_2}{|\xi_1-\xi_2|^{n-2}} $$ $$\le 
  C\sum_{0\le|\alpha_2|\le m-1}\int_{{\bf R}^n}\int_{{\bf R}^n}|\mu(x',\xi_1,\xi_2,y')|^{-2m+2+|\alpha_2|}\left(| x'-\xi_1|^{-n+2+|\alpha_2|}+| y'-\xi_2|^{-n+2+|\alpha_2|}\right)$$ $$\times\left|\partial_{\xi_1+\xi_2}^{\alpha_2}\left(V(h\xi_1)V(h\xi_2)\right)\right|\frac{d\xi_1d\xi_2}{|\xi_1-\xi_2|^{n-2}} $$
  $$\le Ch^{-4}
  \sum_{0\le|\alpha_2|\le m-1}\int_{{\bf R}^n}\int_{{\bf R}^n}|\mu(x,\xi_1,\xi_2,y)|^{-2m+2+|\alpha_2|}\left(|x-\xi_1|^{-n+2+|\alpha_2|}+| y-\xi_2|^{-n+2+|\alpha_2|}\right)$$ $$\times\left|\partial_{\xi_1+\xi_2}^{\alpha_2}(V(\xi_1)V(\xi_2))\right|\frac{d\xi_1d\xi_2}{|\xi_1-\xi_2|^{n-2}} $$
   $$\le Ch^{-4}
  \sum_{p=0}^{m-1}\sum_{p_1+p_2=p}\|V\|_{{\cal C}_{\frac{n+1}{2}+\epsilon'}^{p_1}}\|V\|_{{\cal C}_{\frac{n+1}{2}+\epsilon'}^{p_2}}\int_{{\bf R}^n}\int_{{\bf R}^n}|\mu(x,\xi_1,\xi_2,y)|^{-2m+2+p}$$ $$\times \left(|x-\xi_1|^{-n+2+p}+| y-\xi_2|^{-n+2+p}\right)\langle\xi_1\rangle^{-\frac{n+1}{2}-\epsilon'}\langle\xi_2\rangle^{-\frac{n+1}{2}-\epsilon'}\frac{d\xi_1d\xi_2}{|\xi_1-\xi_2|^{n-2}}. $$
  We need now the following bound proved in the appendix.
  
  \begin{lemma} Let $0\le\ell_1< n-1$, $\ell_1\le\ell_2,\ell_3<n-1$, $\ell_2+\ell_4>n$, $\ell_3+\ell_5>n$. Then we have the bound
  $$\sup_{x,y\in{\bf R}^n}\int_{{\bf R}^n}\int_{{\bf R}^n}|\mu(x,\xi_1,\xi_2,y)|^{-\ell_1}|x-\xi_1|^{-\ell_2}|\xi_1-\xi_2|^{-\ell_3}\langle\xi_1\rangle^{-\ell_4}\langle\xi_2\rangle^{-\ell_5}d\xi_1d\xi_2<+\infty.\eqno{(5.9)}$$
  \end{lemma}
  
   Thus we conclude that
  $$\left\|{\cal B}^{(3)}_1(t,h)\right\|_{L^1\to L^\infty}\le Ch^{m-n}t^{-m}\sum_{p=0}^{m-1}\|V\|_{{\cal C}_{\frac{n+1}{2}+\epsilon'}^{p}}\|V\|_{{\cal C}_{\frac{n+1}{2}+\epsilon'}^{m-1-p}}.\eqno{(5.10)}$$
  If $n$ is odd the desired estimate follows from (5.10) applied with $m=\frac{n-1}{2}$. Let now $n$ be even.
  Then there exists a family of functions, $V_\theta$, $0<\theta\le 1$, such that
  $$\left\|V_\theta\right\|_{{\cal C}_{\frac{n+1}{2}+\epsilon'}^{\frac{n-4}{2}}}+\theta^{-1/2}\left\|V-V_\theta\right\|_{{\cal C}_{\frac{n+1}{2}+\epsilon'}^{\frac{n-4}{2}}}+\theta^{1/2}\left\|V_\theta\right\|_{{\cal C}_{\frac{n+1}{2}+\epsilon'}^{\frac{n-2}{2}}}\le C.\eqno{(5.11)}$$
 It follows from (5.11) that the functions
$$W^{(0)}(\xi_1,\xi_2)=(V(\xi_1)-V_{\theta}(\xi_1))V(\xi_2)+V_{\theta}(\xi_1)(V(\xi_2)-V_{\theta}(\xi_2)),$$
$$W^{(1)}(\xi_1,\xi_2)=V_{\theta}(\xi_1)V_{\theta}(\xi_2)$$
satisfy the bound
$$\left|\langle\xi_1\rangle^{\frac{n+1}{2}+\epsilon'}\langle\xi_2\rangle^{\frac{n+1}{2}+\epsilon'}\partial_{\xi_1+\xi_2}^\alpha W^{(j)}(\xi_1,\xi_2)\right|\le C\theta^{1/2-j},\eqno{(5.12)}$$
for $|\alpha|\le \frac{n-4}{2}+j$. 
 Thus we decompose the function (5.5) as $E^{(0)}+E^{(1)}$, where $E^{(j)}$, $j=0,1$, is defined by replacing in (5.5) the function $V(\xi_1)V(\xi_2)$ by $W^{(j)}(\xi_1,\xi_2)$. Clearly, in this case the desired estimate follows from the following estimates applied with $\theta=h$.

 \begin{lemma} The functions $E^{(j)}$, $j=0,1$, satisfy the estimate
  $$\left|E^{(j)}\right|\le Ch^{-\frac{n+1}{2}}t^{-\frac{n-1}{2}}\left(\frac{\theta}{h}\right)^{1/2-j}.\eqno{(5.13)}$$
  \end{lemma}
  
  {\it Proof.} We first integrate by parts $\frac{n-4}{2}+j$ times with respect to $\xi_1+\xi_2$. Thus we get that $E^{(j)}$ is a linear combination of functions of the form
  $$h^{-n+4}\int_{{\bf R}^n}\int_{{\bf R}^n}e_{\alpha_1,\alpha_2}^{(j)}(t',x',y',\xi_1,\xi_2)q_\alpha^{(\frac{n-4}{2}+j)}(x',\xi_1,\xi_2,y')\partial_{\xi_1+\xi_2}^\beta\left( W^{(j)}(h\xi_1,h\xi_2)\right)d\xi_1d\xi_2,$$
  where
  $$e_{\alpha_1,\alpha_2}^{(j)}=\int e^{i\lambda(t'-|x'-\xi_1|-|y'-\xi_2|)}\lambda^{-\frac{n-4}{2}-j}
  \widetilde\varphi(\lambda)$$ $$\times a^-_{\nu,0}(\lambda|\xi_1-\xi_2|)\partial_{\xi_1}^{\alpha_1}\left(b^-_{\nu,0}(\lambda|x'-\xi_1|)\right)\partial_{\xi_2}^{\alpha_2}\left(b^-_{\nu,0}(\lambda|y'-\xi_2|)\right)d\lambda,$$
  $|\alpha_1|+|\alpha_2|+|\beta|=|\alpha|\le \frac{n-4}{2}+j$. Let $0\le m\le \frac{n}{2}$ be an integer. 
  We now integrate by parts $m$ times with respect to $\lambda$ to obtain
  $$\left|e_{\alpha_1,\alpha_2}^{(j)}\right|\le C(t')^{-m}\sum_{j_1+j_2+j_3\le m}\int_{{\rm supp}\,\widetilde\varphi}\left|a^-_{\nu,j_2}(\lambda|\xi_1-\xi_2|)\right|$$ $$\times 
  \left|\partial_{\xi_1}^{\alpha_1}\left(b^-_{\nu,j_1}(\lambda|x'-\xi_1|)\right)\right|
  \left|\partial_{\xi_2}^{\alpha_2}\left(b^-_{\nu,j_2}(\lambda|y'-\xi_2|)\right)\right|d\lambda$$
   $$\le C(t')^{-m}|\xi_1-\xi_2|^{-n+2}\sum_{j_1+j_2+j_3\le m}\langle x'-\xi_1\rangle^{j_1-\frac{n-1}{2}-|\alpha_1|}\langle y'-\xi_2\rangle^{j_3-\frac{n-1}{2}-|\alpha_2|}$$
  $$\le C(t')^{-m}|\xi_1-\xi_2|^{-n+2}\left(\langle x'-\xi_1\rangle+\langle y'-\xi_2\rangle\right)^m\langle x'-\xi_1\rangle^{-\frac{n-1}{2}-|\alpha_1|}\langle y'-\xi_2\rangle^{-\frac{n-1}{2}-|\alpha_2|},
  \eqno{(5.14)}$$
  where we have used (4.11) and (5.6). By interpolation, (5.14) holds for all real $0\le m\le \frac{n}{2}$ and in particular for $m=\frac{n-1}{2}$. Hence, using this together with (5.4), we obtain
  $$\left|e_{\alpha_1,\alpha_2}^{(j)}q_\alpha^{(\frac{n-4}{2}+j)}\right|$$ $$\le C(t')^{-\frac{n-1}{2}}|\xi_1-\xi_2|^{-n+2}\left(|x'-\xi_1|^{-1}+| y'-\xi_2|^{-1}\right)^{\frac{n-1}{2}+|\alpha_1|+|\alpha_2|}\left|q_\alpha^{(\frac{n-4}{2}+j)}(x',\xi_1,\xi_2,y')
  \right|$$ $$\le C(t')^{-\frac{n-1}{2}}|\xi_1-\xi_2|^{-n+2}|\mu|^{-n+4-2j+|\beta|}\left(|x'-\xi_1|^{-1}+| y'-\xi_2|^{-1}\right)^{n-\frac{5}{2}+j-|\beta|}.\eqno{(5.15)}$$
  By (5.12) and (5.15), we obtain
  $$(t/h)^{\frac{n-1}{2}}h^{n}\left|E^{(j)}\right|\le C(\theta/h)^{1/2-j}\sum_{p=0}^{\frac{n-4}{2}+j}\int_{{\bf R}^n}\int_{{\bf R}^n}|\mu(x,\xi_1,\xi_2,y)|^{-n+4-2j+p}$$ $$\times \left(|x-\xi_1|^{-n+5/2-j+p}+| y-\xi_2|^{-n+5/2-j+p}\right)\langle\xi_1\rangle^{-\frac{n+1}{2}-\epsilon'}\langle\xi_2\rangle^{-\frac{n+1}{2}-\epsilon'}\frac{d\xi_1d\xi_2}{|\xi_1-\xi_2|^{n-2}},$$
  which together with Lemma 5.3 imply (5.13).
  \eproof

  \begin{prop} Let $V$ satisfy (1.6). Then in dimensions $n=6,7$ the operator ${\cal B}^{(1)}$ satisfies the estimate
$$\left\|{\cal B}^{(1)}(t,h)\right\|_{L^1\to L^\infty}\le C_\epsilon h^{-\frac{n+1}{2}-\epsilon}t^{-\frac{n-1}{2}},\eqno{(5.16)}$$
for every $0<\epsilon\ll 1$.
\end{prop}

{\it Proof.} Making a change of variables $\xi_1\to h\xi_1$, $\xi_2\to h\xi_2$, we can write the kernel of the operator ${\cal B}^{(1)}(t,h)$ in the form
$$h^{-n+4}\int_{{\bf R}^n}\int_{{\bf R}^n}B_1^{(1)}(|x'-\xi_1|,|\xi_1-\xi_2|,|\xi_2-y'|,t')V(h\xi_1)V(h\xi_2)d\xi_1d\xi_2,\eqno{(5.17)}$$
where $x'=x/h$, $y'=y/h$, $t'=t/h$.  Let $1\le m\le\frac{n-1}{2}$ be an integer. Integrating by parts $m$ times with respect to the variable $\lambda$ we can write the function $B_1^{(1)}$ in the form
$$\sum_{j_1+j_2+j_3\le m}B^{(1)}_{1,j_1,j_2,j_3},$$
where 
$$B^{(1)}_{1,j_1,j_2,j_3}(\sigma_1,\sigma_2,\sigma_3,t)=t^{-m}\int e^{i(t-(\sigma_1+\sigma_2+\sigma_3))\lambda}\varphi_{j_1,j_2,j_3}(\lambda)b_{\nu,j_1}^-(\sigma_1\lambda)b_{\nu,j_2}^-(\sigma_2\lambda)b_{\nu,j_3}^-(\sigma_3\lambda)d\lambda,$$
$\varphi_{j_1,j_2,j_3}\in C_0^\infty((0,+\infty))$, $b_{\nu,j}^-(z)$ being defined in the previous section. We have to show that the operator with kernel
 $$h^{-n+4}\int_{{\bf R}^n}\int_{{\bf R}^n}B_{1,j_1,j_2,j_3}^{(1)}(|x'-\xi_1|,|\xi_1-\xi_2|,|\xi_2-y'|,t')V(h\xi_1)V(h\xi_2)d\xi_1d\xi_2\eqno{(5.18)}$$
 satisfies (5.16). Set
$$\Phi(x,\xi_1,\xi_2,y)=|x-\xi_1|+|\xi_1-\xi_2|+|y-\xi_2|,$$ 
 $$\rho_1:= \nabla_{\xi_1} \Phi(x,\xi_1,\xi_2,y)=\rho(x,\xi_1,\xi_2),\quad \rho_2:=\nabla_{\xi_2} \Phi(x,\xi_1,\xi_2,y)=\rho(\xi_1,\xi_2,y),$$
$$ w_1=|x-\xi_1|^{-1}+|\xi_1-\xi_2|^{-1},\quad w_2=|\xi_1-\xi_2|^{-1}+|y-\xi_2|^{-1}.$$
We will need the following bounds proved in the appendix.
  
  \begin{lemma} Let $0\le\ell_j< n-1$, $j=1,2$, $0\le\ell_j<n$,  $\ell_j+\ell_5>n$, $j=3,4$. If either
  $$\min\{\ell_1,\ell_2\}\le \min\{\ell_3,\ell_4\},\quad \max\{\ell_1,\ell_2\}\le \max\{\ell_3,\ell_4\},\quad \ell_2\le\ell_4\, ({\rm if}\, \ell_3>\ell_4),\eqno{(5.19)}$$
  or
   $$\ell_1\le\frac{\ell_3+\ell_4}{2},\quad\ell_2\le\min\left\{\frac{\ell_3+\ell_4}{2},\ell_4\right\},\eqno{(5.20)}$$
 we have the bound
  $$\sup_{x,y\in{\bf R}^n}\int_{{\bf R}^n}\int_{{\bf R}^n}|\rho(x,\xi_1,\xi_2)|^{-\ell_1}|\rho(\xi_1,\xi_2,y)|^{-\ell_2}|x-\xi_1|^{-\ell_3}|\xi_1-\xi_2|^{-\ell_4}\langle\xi_1\rangle^{-\ell_5}\langle\xi_2\rangle^{-\ell_5}d\xi_1d\xi_2<+\infty.\eqno{(5.21)}$$
  If either
  $$\min\{\ell_1,\ell_2\}\le \min\{\ell_3,\ell_4\},\quad \max\{\ell_1,\ell_2\}\le \max\{\ell_3,\ell_4\},\quad \ell_1\le\ell_3\, ({\rm if}\, \ell_4>\ell_3),\eqno{(5.22)}$$
  or
  $$\ell_1\le\min\left\{\frac{\ell_3+\ell_4}{2},\ell_3\right\},\quad\ell_2\le\frac{\ell_3+\ell_4}{2},\eqno{(5.23)}$$
 we have the bound
  $$\sup_{x,y\in{\bf R}^n}\int_{{\bf R}^n}\int_{{\bf R}^n}|\rho(x,\xi_1,\xi_2)|^{-\ell_1}|\rho(\xi_1,\xi_2,y)|^{-\ell_2}|\xi_1-\xi_2|^{-\ell_3}|y-\xi_2|^{-\ell_4}\langle\xi_1\rangle^{-\ell_5}\langle\xi_2\rangle^{-\ell_5}d\xi_1d\xi_2<+\infty.\eqno{(5.24)}$$
  If $$\ell_1\le\ell_3,\quad\ell_2\le\ell_4,\eqno{(5.25)}$$
 we have the bound
  $$\sup_{x,y\in{\bf R}^n}\int_{{\bf R}^n}\int_{{\bf R}^n}|\rho(x,\xi_1,\xi_2)|^{-\ell_1}|\rho(\xi_1,\xi_2,y)|^{-\ell_2}|x-\xi_1|^{-\ell_3}|y-\xi_2|^{-\ell_4}\langle\xi_1\rangle^{-\ell_5}\langle\xi_2\rangle^{-\ell_5}d\xi_1d\xi_2<+\infty.\eqno{(5.26)}$$
  \end{lemma}
 Our task is to gain an additional factor $O(h^{n-3})$. This can be achieved by integrating by parts with respect to $\xi_1$ and $\xi_2$. This procedure, however, leads to integrals with singularities which could be a priori too strong. Our aim is to perform the integration by parts in such a way that at the end we get singular integrals covered by Lemma 5.6. This is far from being obvious and probably impossible to do in high dimensions. In our case, however, this is relatively easy if we allow an $\epsilon$ loss in $h$. 
We will first consider the case $n=7$. Take $m=\frac{n-1}{2}=3$. We are going to integrate by parts once with respect to $\xi_2$ and then twice with respect to $\xi_1$ using the identity
$$\lambda^3\left(\Lambda_{\xi_1}^*\right)^2\Lambda_{\xi_2}^*=\sum_{0\le|\alpha_1|\le 2}r^{(2)}_{\alpha_1}(x,\xi_1,\xi_2)\partial_{\xi_1}^{\alpha_1}\sum_{0\le|\alpha_2|\le 1}r^{(1)}_{\alpha_2}(\xi_1,\xi_2,y)\partial_{\xi_2}^{\alpha_2}$$ 
$$=\sum_{|\alpha_1|\le 2,\,|\alpha_2|\le 1}g_{\alpha_1,\alpha_2}(x,\xi_1,\xi_2,y)\partial_{\xi_1}^{\alpha_1}\partial_{\xi_2}^{\alpha_2}$$ $$=\sum_{|\alpha_1|\le 2,\,|\alpha_2|\le 1}\sum_{\kappa_1,\kappa_2\in\Theta(\alpha_1,\alpha_2)}g_{\alpha_1,\alpha_2}^{\kappa_1,\kappa_2}(x,\xi_1,\xi_2,y)\partial_{\xi_1}^{\alpha_1}\partial_{\xi_2}^{\alpha_2}\eqno{(5.27)}$$
where $g_{\alpha_1,\alpha_2}^{\kappa_1,\kappa_2}=O_{w_1,w_2}\left(|\rho_1|^{-2-\kappa_1}|\rho_2|^{-1-\kappa_2}\right)$, $\Theta(\alpha_1,\alpha_2)$ denotes the set of all integers $\kappa_1,\kappa_2\ge 0$ such that $\kappa_1+\kappa_2=3-|\alpha_1|-|\alpha_2|$, $\kappa_1\le 2-|\alpha_1|$. More precisely, using (4.9) one can check that the functions $g_{\alpha_1,\alpha_2}^{\kappa_1,\kappa_2}$ satisfy the bounds (with $0\le|\beta_1|,|\beta_2|\le 1$):
$$\left|\partial_{\xi_1}^{\beta_1} \partial_{\xi_2}^{\beta_2}g_{\alpha_1,\alpha_2}^{\kappa_1,\kappa_2}\right|\le C |\rho_1|^{-2-\kappa_1}|\rho_2|^{-1-\kappa_2}w_1^{\kappa_1}w_2^{\kappa_2}\left(w_1|\rho_1|^{-1}+w_2|\rho_2|^{-1}\right)^{|\beta_1|+|\beta_2|}.\eqno{(5.28)}$$
Thus, using (5.27) we obtain that (5.18) is a linear combination of functions of the form
$$t^{-3}\int\int_{{\bf R}^7}\int_{{\bf R}^7}e^{i\lambda(t'-|x'-\xi_1|-|\xi_1-\xi_2|-|y'-\xi_2|)}
 \lambda^{-3}\varphi_{j_1,j_2,j_3}(\lambda)g_{\alpha_1,\alpha_2}^{\kappa_1,\kappa_2}(x',\xi_1,\xi_2,y')$$ 
$$\times\partial_{\xi_1}^{\beta_1}\left(b_{\nu,j_1}^-(\lambda|x'-\xi_1|)\right)\partial_{\xi_1-\xi_2}^{\beta_2}\left(b^-_{\nu,j_2}(\lambda|\xi_1-\xi_2|)\right)\partial_{\xi_2}^{\beta_3}\left(b_{\nu,j_3}^-(\lambda|y'-\xi_2|)\right)$$ $$\times\partial_{\xi_1}^{\gamma_1}\left(V(h\xi_1)\right)\partial_{\xi_2}^{\gamma_2}\left(V(h\xi_2)\right)d\xi_1d\xi_2 d\lambda,\eqno{(5.29)}$$
where $0\le|\alpha_1|\le 2$, $0\le|\alpha_2|\le 1$, $\kappa_1,\kappa_2\in\Theta(\alpha_1,\alpha_2)$,  $|\beta_1|+|\beta_2|+|\beta_3|+|\gamma_1|+|\gamma_2|=|\alpha_1|+|\alpha_2|$. 
Clearly,
we have $0\le|\gamma_1|\le 2$, $0\le|\gamma_2|\le 1$. Set
$$A_{k_1,k_2}(x,\xi_1,\xi_2,y)=|x-\xi_1|^{-k_1}\left(|\xi_1-\xi_2|^{-k_2}+|y-\xi_2|^{-k_2}\right)$$ $$ +|\xi_1-\xi_2|^{-k_1}\left(|x-\xi_1|^{-k_2}+|y-\xi_2|^{-k_2}\right)+|y-\xi_2|^{-k_1}\left(|x-\xi_1|^{-k_2}+|\xi_1-\xi_2|^{-k_2}\right).$$
Note that there exists a constant $c>0$ such that $$|x'-\xi_1|\ge c, \quad |\xi_1-\xi_2|\ge c, \quad|y'-\xi_2|\ge c,$$ as long as $\lambda|x'-\xi_1|\in$ supp$\,b^-_{\nu,j_1}$, $\lambda|\xi_1-\xi_2|\in$ supp$\,b^-_{\nu,j_2}$, $\lambda|y'-\xi_2|\in$ supp$\,b^-_{\nu,j_3}$ and 
$\lambda\in$ supp$\,\varphi_{j_1,j_2,j_3}$. Using this together with (4.11) we obtain the bound (in any dimension $n$)
$$\left|\partial_{\xi_1}^{\beta_1}\left(b_{\nu,j_1}^-(\lambda|x'-\xi_1|)\right)\partial_{\xi_1-\xi_2}^{\beta_2}\left(b^-_{\nu,j_2}(\lambda|\xi_1-\xi_2|)\right)\partial_{\xi_2}^{\beta_3}\left(b_{\nu,j_3}^-(\lambda|y'-\xi_2|)\right)\right|$$
$$\le C|x'-\xi_1|^{j_1-\frac{n-1}{2}-|\beta_1|}|\xi_1-\xi_2|^{j_2-\frac{n-1}{2}-|\beta_2|}|y'-\xi_2|^{j_3-\frac{n-1}{2}-|\beta_3|}$$
$$\le \frac{C\left(|x'-\xi_1|+|\xi_1-\xi_2|+|y'-\xi_2|\right)^m}{|x'-\xi_1|^{\frac{n-1}{2}}|\xi_1-\xi_2|^{\frac{n-1}{2}}|y'-\xi_2|^{\frac{n-1}{2}}}
\left(|x'-\xi_1|^{-1}+|\xi_1-\xi_2|^{-1}+|y'-\xi_2|^{-1}\right)^p$$
 $$\le \frac{C\left(|x'-\xi_1|^{\frac{n-1}{2}}+|\xi_1-\xi_2|^{\frac{n-1}{2}}+|y'-\xi_2|^{\frac{n-1}{2}}\right)}{|x'-\xi_1|^{\frac{n-1}{2}}|\xi_1-\xi_2|^{\frac{n-1}{2}}|y'-\xi_2|^{\frac{n-1}{2}}}
\left(|x'-\xi_1|^{-p}+|\xi_1-\xi_2|^{-p}+|y'-\xi_2|^{-p}\right)$$
$$\le CA_{\frac{n-1}{2}+p,\frac{n-1}{2}}(x',\xi_1,\xi_2,y')\le CA_{\frac{n-1}{2}+p(1-\epsilon),\frac{n-1}{2}}(x',\xi_1,\xi_2,y'),\eqno{(5.30)}$$
for every $0<\epsilon\ll 1$, where we have put $p=|\beta_1|+|\beta_2|+|\beta_3|$. We now proceed as follows. If $|\gamma_1|=2$ we integrate by parts once with respect to $\xi_2$ in the integral (5.29). When  $|\gamma_1|\le 1$, we integrate by parts once with respect to $\xi_2$ in the integral (5.29) if $1+\kappa_1\ge\kappa_2$, and with respect to $\xi_1$ if $1+\kappa_1<\kappa_2$. This procedure together with (4.9), (5.28) and (5.30) (with $\frac{n-1}{2}=3$) lead to the bound 
$$\left|(5.29)\right|\le 
 Ct^{-3}\sum_{j=0}^1\sum_{0\le|\gamma_1|,|\gamma_2|\le 2}\sum_{p=0}^{4-|\gamma_1|-|\gamma_2|}\sum_{\kappa_1,\kappa_2\in\Theta_j^\sharp(p)}\int_{{\bf R}^7}\int_{{\bf R}^7} |\rho_1|^{-2-\kappa_1-j}|\rho_2|^{-2-\kappa_2+j}$$ $$\times \left(|x'-\xi_1|^{-\kappa_1(1-\epsilon)}+|\xi_1-\xi_2|^{-\kappa_1(1-\epsilon)}\right)\left(|\xi_1-\xi_2|^{-\kappa_2(1-\epsilon)}+|y'-\xi_2|^{-\kappa_2(1-\epsilon)}\right)$$
$$\times A_{3+p(1-\epsilon) ,3}(x',\xi_1,\xi_2,y')\left|\partial_{\xi_1}^{\gamma_1}\left(V(h\xi_1)\right)\right|\left|\partial_{\xi_2}^{\gamma_2}\left(V(h\xi_2)\right)\right|d\xi_1d\xi_2$$
 $$\le Ch^{-4-4\epsilon}t^{-3}\sum_{j=0}^1\sum_{0\le|\gamma_1|,|\gamma_2|\le 2}\sum_{p=0}^{4-|\gamma_1|-|\gamma_2|}\sum_{\kappa_1,\kappa_2\in\Theta_j^\sharp(p)}\int_{{\bf R}^7}\int_{{\bf R}^7} |\rho_1|^{-2-\kappa_1-j}|\rho_2|^{-2-\kappa_2+j}$$ $$\times \left(|x-\xi_1|^{-\kappa_1(1-\epsilon)}+|\xi_1-\xi_2|^{-\kappa_1(1-\epsilon)}\right)\left(|\xi_1-\xi_2|^{-\kappa_2(1-\epsilon)}+|y-\xi_2|^{-\kappa_2(1-\epsilon)}\right)$$
$$\times A_{3+p(1-\epsilon),3}(x,\xi_1,\xi_2,y)\left|\partial_{\xi_1}^{\gamma_1}V(\xi_1)\right|\left|\partial_{\xi_2}^{\gamma_2}V(\xi_2)\right|d\xi_1d\xi_2$$
 $$\le Ch^{-4-4\epsilon}t^{-3}\sum_{j=0}^1\sum_{p=0}^{4}\sum_{\kappa_1,\kappa_2\in\Theta_j^\sharp(p)}\int_{{\bf R}^7}\int_{{\bf R}^7} |\rho_1|^{-2-\kappa_1-j}|\rho_2|^{-2-\kappa_2+j}$$ $$\times \left(|x-\xi_1|^{-\kappa_1(1-\epsilon)}+|\xi_1-\xi_2|^{-\kappa_1(1-\epsilon)}\right) \left(|\xi_1-\xi_2|^{-\kappa_2(1-\epsilon)}+|y-\xi_2|^{-\kappa_2(1-\epsilon)}\right)$$
$$\times A_{3+p(1-\epsilon),3}(x,\xi_1,\xi_2,y)\langle\xi_1\rangle^{-4-\epsilon'}\langle\xi_2\rangle^{-4-\epsilon'}d\xi_1d\xi_2\eqno{(5.31)}$$
for every $0<\epsilon\ll 1$ and some $\epsilon'>0$, where $\Theta_0^\sharp(p)$ denotes the set of all integers $0\le\kappa_1,\kappa_2\le 3$ such that $\kappa_1+\kappa_2\le 4-p$, $\Theta_1^\sharp(p)$ denotes the set of all integers $0\le\kappa_1\le 1$, $0\le\kappa_2\le 4$ such that $\kappa_1+\kappa_2\le 4-p$. Note that the $\epsilon$ loss in (5.31) allows to avoid non-integrable singularities. 
It is not hard to see that the right-hand side of (5.31) is bounded by singular integrals satisfying the conditions of Lemma 5.6, which yields the desired estimate in this case.

 Let now $n=6$. Then there exists a family of functions, $V_\theta$, $0<\theta\le 1$, such that
  $$\left\|V_\theta\right\|_{{\cal C}_{7/2+\epsilon'}^1}+\theta^{-1/2}\left\|V-V_\theta\right\|_{{\cal C}_{7/2+\epsilon'}^1}+\theta^{1/2}\left\|V_\theta\right\|_{{\cal C}_{7/2+\epsilon'}^2}\le C,\eqno{(5.32)}$$
with some constants $C,\epsilon'>0$ independent of $\theta$. In this case we would like to gain a factor $O(h^3)$. We will modify a little bit the integration by parts scheme used in the case $n=7$ above. Note first that, as above, we have
$$\lambda^2\Lambda_{\xi_1}^*\Lambda_{\xi_2}^*=\sum_{0\le|\alpha_1|\le 1}r^{(1)}_{\alpha_1}(x,\xi_1,\xi_2)\partial_{\xi_1}^{\alpha_1}\sum_{0\le|\alpha_2|\le 1}r^{(1)}_{\alpha_2}(\xi_1,\xi_2,y)\partial_{\xi_2}^{\alpha_2}$$
$$=\sum_{0\le|\alpha_1|,|\alpha_2|\le 1}f_{\alpha_1,\alpha_2}(x,\xi_1,\xi_2,y)\partial_{\xi_1}^{\alpha_1}\partial_{\xi_2}^{\alpha_2}$$
$$=\sum_{|\alpha_1|,|\alpha_2|\le 1}\sum_{\kappa_1,\kappa_2\in\Omega(\alpha_1,\alpha_2)}f_{\alpha_1,\alpha_2}^{\kappa_1,\kappa_2}(x,\xi_1,\xi_2,y)\partial_{\xi_1}^{\alpha_1}\partial_{\xi_2}^{\alpha_2},\eqno{(5.33)}$$
with functions $f_{\alpha_1,\alpha_2}$ satisfying the bounds (with $0\le|\beta_1|,|\beta_2|\le 2$):
$$\left|\partial_{\xi_1}^{\beta_1} \partial_{\xi_2}^{\beta_2}f_{\alpha_1,\alpha_2}^{\kappa_1,\kappa_2}\right|\le C |\rho_1|^{-1-\kappa_1}|\rho_2|^{-1-\kappa_2}w_1^{\kappa_1}w_2^{\kappa_2}\left(w_1|\rho_1|^{-1}+w_2|\rho_2|^{-1}\right)^{|\beta_1|+|\beta_2|},\eqno{(5.34)}$$
where $\Omega(\alpha_1,\alpha_2)$ denotes the set of all integers $\kappa_1$, $\kappa_2$ satisfying $\kappa_1+\kappa_2=2-|\alpha_1|-|\alpha_2|$, $0\le \kappa_1\le 1-|\alpha_1|$.
We integrate by parts in (5.17) successively with respect to $\xi_2$ and $\xi_1$ using (5.33). Thus we get that the function (5.17) is a linear combination of functions of the form
$$h^{-2}\int_{{\bf R}^6}\int_{{\bf R}^6}u_{\beta_1,\beta_2,\beta_3}(t',x',y',\xi_1,\xi_2)f_{\alpha_1,\alpha_2}^{\kappa_1,\kappa_2}(x',\xi_1,\xi_2,y')\partial_{\xi_1}^{\gamma_1}\left( V(h\xi_1)\right)\partial_{\xi_2}^{\gamma_2}\left(V(h\xi_2)\right)d\xi_1d\xi_2,\eqno{(5.35)}$$
where
 $$u_{\beta_1,\beta_2,\beta_3}=\int e^{i\lambda(t'-|x'-\xi_1|-|\xi_1-\xi_2|-|y'-\xi_2|)}
 \lambda^{-2}\widetilde\varphi(\lambda)$$ $$\times\partial_{\xi_1}^{\beta_1}\left( b_{\nu}^-(\lambda|x'-\xi_1|)\right)\partial_{\xi_1-\xi_2}^{\beta_2}\left(b^-_{\nu}(\lambda|\xi_1-\xi_2|)\right)\partial_{\xi_2}^{\beta_3}\left(b_{\nu}^-(\lambda|y'-\xi_2|)\right) d\lambda,$$
where $0\le|\alpha_1|,|\alpha_2|\le 1$, $\kappa_1,\kappa_2\in\Omega(\alpha_1,\alpha_2)$, $|\beta_1|+|\beta_2|+|\beta_3|+|\gamma_1|+|\gamma_2|=|\alpha_1|+|\alpha_2|$. Clearly, $|\gamma_1|\le|\alpha_1|$, $|\gamma_2|\le|\alpha_2|$. We will consider two cases.

Case 1. $|\gamma_1|=|\gamma_2|=1$. Then we have $|\alpha_1|=|\alpha_2|=1$, $\beta_1=\beta_2=\beta_3=0$, $\kappa_1=\kappa_2=0$. Set
$$Y_{k_1,k_2}=h^{-2}\int_{{\bf R}^6}\int_{{\bf R}^6}u_{0,0,0}(t',x',y',\xi_1,\xi_2)f_{\alpha_1,\alpha_2}^{0,0}(x',\xi_1,\xi_2,y')$$ 
$$\times\partial_{\xi_1}^{\gamma_1}\left(V(h\xi_1)\right)\partial_{\xi_2}^{\gamma_2}\left(V(h\xi_2)\right)\phi_{k_1}\left(|\rho_1|^{-1}\right)\phi_{k_2}\left(|\rho_2|^{-1}\right)d\xi_1d\xi_2,$$
with $\phi_k$ satisfying (4.17). Clearly, it suffices to show that
$$\left|Y_{k_1,k_2}\right|\le C 2^{-\epsilon_0(k_1+k_2)}h^{-7/2}t^{-5/2},\eqno{(5.36)}$$
for some $\epsilon_0>0$. To do so, we decompose the function $\partial_{\xi_1}^{\gamma_1}V(\xi_1)\partial_{\xi_2}^{\gamma_2}V(\xi_2)$ as $\sum_{j=1}^4W_j$, where
$$W_1(\xi_1,\xi_2)=\partial_{\xi_1}^{\gamma_1}\left(V(\xi_1)-V_\theta(\xi_1)\right)\partial_{\xi_2}^{\gamma_2}\left(V(\xi_2)-V_\theta(\xi_2)\right),$$
 $$W_2(\xi_1,\xi_2)=\partial_{\xi_1}^{\gamma_1}\left(V(\xi_1)-V_\theta(\xi_1)\right)\partial_{\xi_2}^{\gamma_2}V_\theta(\xi_2),$$
 $$W_3(\xi_1,\xi_2)=\partial_{\xi_1}^{\gamma_1}V_\theta(\xi_1)\partial_{\xi_2}^{\gamma_2}\left(V(\xi_2)-V_\theta(\xi_2)\right),$$
 $$W_4(\xi_1,\xi_2)=\partial_{\xi_1}^{\gamma_1}V_\theta(\xi_1)\partial_{\xi_2}^{\gamma_2}V_\theta(\xi_2).$$
 It follows from (5.32) that
$$\theta^{-1}\left|\langle\xi_1\rangle^{7/2+\epsilon'}\langle\xi_2\rangle^{7/2+\epsilon'}W_1(\xi_1,\xi_2)\right|+\left|\langle\xi_1\rangle^{7/2+\epsilon'}\langle\xi_2\rangle^{7/2+\epsilon'}\partial_{\xi_2}^{\alpha_2}W_2(\xi_1,\xi_2)\right|$$ $$+\left|\langle\xi_1\rangle^{7/2+\epsilon'}\langle\xi_2\rangle^{7/2+\epsilon'}\partial_{\xi_1}^{\alpha_1}W_3(\xi_1,\xi_2)\right|+\theta\left|\langle\xi_1\rangle^{7/2+\epsilon'}\langle\xi_2\rangle^{7/2+\epsilon'}\partial_{\xi_1}^{\alpha_1}\partial_{\xi_2}^{\alpha_2}W_4(\xi_1,\xi_2)\right|\le C,\eqno{(5.37)}$$
for all $0\le|\alpha_1|,|\alpha_2|\le 1$. Write
  $$Y_{k_1,k_2}=\sum_{j=1}^4Y^{(j)}_{k_1,k_2},$$ where
$Y^{(j)}_{k_1,k_2}$ is defined by replacing in the definition of $Y_{k_1,k_2}$ the function $\partial_{\xi_1}^{\gamma_1}(V(h\xi_1))\partial_{\xi_2}^{\gamma_2}(V(h\xi_2))$ by $h^2W_j(h\xi_1,h\xi_2)$. 

\begin{lemma} The functions $Y_{k_1,k_2}^{(j)}$, $j=1,2,3,4$, satisfy the estimates
  $$\left|Y_{k_1,k_2}^{(1)}\right|\le C2^{-\epsilon_0(k_1+k_2)} h^{-7/2}t^{-5/2}\left(\frac{\theta }{h2^{k_1+k_2}}\right),\eqno{(5.38)}$$
   $$\left|Y_{k_1,k_2}^{(2)}\right|\le C2^{-\epsilon_0(k_1+k_2)} h^{-7/2}t^{-5/2},\eqno{(5.39)}$$
   $$\left|Y_{k_1,k_2}^{(3)}\right|\le C2^{-\epsilon_0(k_1+k_2)} h^{-7/2}t^{-5/2},\eqno{(5.40)}$$
   $$\left|Y_{k_1,k_2}^{(4)}\right|\le C2^{-\epsilon_0(k_1+k_2)} h^{-7/2}t^{-5/2}\left(\frac{h2^{k_1+k_2}}{\theta}\right),\eqno{(5.41)}$$
   for some $\epsilon_0>0$.
  \end{lemma}
  
  {\it Proof.} Let $0\le m\le 3$ be an integer. We integrate by parts $m$ times with respect to $\lambda$ to obtain
  $$\left|u_{0,0,0}\right|\le C(t')^{-m}\sum_{j_1+j_2+j_3\le m}\int_{{\rm supp}\,\widetilde\varphi} 
  \left|b^-_{\nu,j_1}(\lambda|x'-\xi_1|)\right|\left|b^-_{\nu,j_2}(\lambda|\xi_1-\xi_2|)\right|
  \left|b^-_{\nu,j_2}(\lambda|y'-\xi_2|)\right|d\lambda$$
   $$\le C(t')^{-m}\sum_{j_1+j_2+j_3\le m}\langle x'-\xi_1\rangle^{j_1-\frac{5}{2}}\langle\xi_1-\xi_2\rangle^{j_2-\frac{5}{2}}\langle y'-\xi_2\rangle^{j_3-\frac{5}{2}}$$
  $$\le C(t')^{-m}\left(\langle x'-\xi_1\rangle+\langle \xi_1-\xi_2\rangle+\langle y'-\xi_2\rangle\right)^m\langle x'-\xi_1\rangle^{-\frac{5}{2}}\langle \xi_1-\xi_2\rangle^{-\frac{5}{2}}\langle y'-\xi_2\rangle^{-\frac{5}{2}},
  \eqno{(5.42)}$$
  where we have used (4.11). Clearly, (5.42) holds for all real $0\le m\le 3$ and in particular for $m=5/2$. Using this together with (5.34) and (5.37), we get
   $$\left|Y_{k_1,k_2}^{(1)}\right|\le C\theta h^{-9/2}t^{-5/2}2^{-k_1(\epsilon_1+1)}2^{-k_2(\epsilon_2+1)}$$
   $$\times \int_{{\bf R}^6}\int_{{\bf R}^6} |\rho_1|^{-2-\epsilon_1}|\rho_2|^{-2-\epsilon_2}A_{5/2,5/2}(x,\xi_1,\xi_2,y)
   \langle\xi_1\rangle^{-7/2-\epsilon'}\langle\xi_2\rangle^{-7/2-\epsilon'} d\xi_1d\xi_2.$$
 Thus (5.38) follows from this estimate and Lemma 5.6, provided $0<\epsilon_1,\epsilon_2\ll 1$ are properly chosen. To get (5.39) we first integrate by parts once with respect to $\xi_2$. Thus we obtain that $Y_{k_1,k_2}^{(2)}$ is a linear combination of functions of the form
 $$\int_{{\bf R}^6}\int_{{\bf R}^6}u^{(1)}_{\alpha'}(t',x',y',\xi_1,\xi_2)v^{(1)}_{\beta'}(x',y',\xi_1,\xi_2)
\partial_{\xi_2}^{\gamma'}\left(W_2(h\xi_1,h\xi_2)\right)d\xi_1d\xi_2,$$
where
$$u^{(1)}_{\alpha'}=\int e^{i\lambda(t'-|x'-\xi_1|-|\xi_1-\xi_2|-|y'-\xi_2|)}
 \lambda^{-3}\widetilde\varphi(\lambda)\partial_{\xi_2}^{\alpha'}\left( b_{\nu}^-(\lambda|x'-\xi_1|)b^-_{\nu}(\lambda|\xi_1-\xi_2|)b_{\nu}^-(\lambda|y'-\xi_2|)\right) d\lambda,$$
 $$v^{(1)}_{\beta'}=r^{(1)}_{\beta'}(\xi_1,\xi_2,y')\partial_{\xi_2}^{\beta'}\left(f_{\alpha_1,\alpha_2}^{0,0}(x',\xi_1,\xi_2,y')\phi_{k_1}\left(|\rho_1|^{-1}\right)\phi_{k_2}\left(|\rho_2|^{-1}\right)\right),$$
$|\alpha'|+|\beta'|+|\gamma'|=1$. As above, integrating by parts with respect to $\lambda$, we obtain 
$$\left|u^{(1)}_{\alpha'}\right|\le C(t')^{-5/2}A_{5/2+p,5/2}(x',\xi_1,\xi_2,y'),\eqno{(5.43)}$$
where $p=|\alpha'|$. On the other hand, by (4.9) and (5.34), we have
$$\left|v^{(1)}_{\beta'}\right|\le C2^{-k_1\epsilon_1}2^{-k_2\epsilon_2}|\rho_1|^{-1-\epsilon_1}|\rho_2|^{-2-\kappa-\epsilon_2}\left(|\xi_1-\xi_2|^{-\kappa}+|\xi_2-y'|^{-\kappa}\right),\eqno{(5.44)}$$
where $\kappa=|\beta'|$. By (5.37), (5.43) and (5.44), we obtain
$$\left|Y_{k_1,k_2}^{(2)}\right|\le Ch^{-7/2}t^{-5/2}2^{-k_1\epsilon_1}2^{-k_2\epsilon_2} \sum_{p=0}^1\sum_{\kappa=0}^{1-p}\int_{{\bf R}^6}\int_{{\bf R}^6}  |\rho_1|^{-1-\epsilon_1}|\rho_2|^{-2-\kappa-\epsilon_2}$$
   $$\times\left(|\xi_1-\xi_2|^{-\kappa}+|\xi_2-y|^{-\kappa}\right) A_{5/2+p,5/2}(x,\xi_1,\xi_2,y)
   \langle\xi_1\rangle^{-7/2-\epsilon'}\langle\xi_2\rangle^{-7/2-\epsilon'} d\xi_1d\xi_2,$$
 which together with Lemma 5.6 yield (5.39). The estimate (5.40) is proved in the same way switching the roles of $\xi_1$ and $\xi_2$. To get (5.41) we integrate by parts successively with respect to $\xi_2$ and $\xi_1$ using (5.33). Thus we obtain that $Y_{k_1,k_2}^{(4)}$ is a linear combination of functions of the form
 $$\int_{{\bf R}^6}\int_{{\bf R}^6}u^{(2)}_{\alpha'_1,\alpha'_2}(t',x',y',\xi_1,\xi_2)v^{(2)}_{\beta'_1,\beta'_2}(x',y',\xi_1,\xi_2)
\partial_{\xi_1}^{\gamma'_1}\partial_{\xi_2}^{\gamma'_2}\left(W_4(h\xi_1,h\xi_2)\right)d\xi_1d\xi_2,$$
where
$$u^{(2)}_{\alpha'_1,\alpha'_2}=$$ $$\int e^{i\lambda(t'-|x'-\xi_1|-|\xi_1-\xi_2|-|y'-\xi_2|)}
 \lambda^{-4}\widetilde\varphi(\lambda)\partial_{\xi_1}^{\alpha'_1}\partial_{\xi_2}^{\alpha'_2}\left( b_{\nu}^-(\lambda|x'-\xi_1|)b^-_{\nu}(\lambda|\xi_1-\xi_2|)b_{\nu}^-(\lambda|y'-\xi_2|)\right) d\lambda,$$
 $$v^{(2)}_{\beta'_1,\beta'_2}=f_{\beta'_1,\beta'_2}^{\kappa'_1,\kappa'_2}(x',\xi_1,\xi_2,y')\partial_{\xi_1}^{\beta'_1}\partial_{\xi_2}^{\beta'_2}\left(f_{\alpha_1,\alpha_2}^{0,0}(x',\xi_1,\xi_2,y')\phi_{k_1}\left(|\rho_1|^{-1}\right)\phi_{k_2}\left(|\rho_2|^{-1}\right)\right),$$
$|\alpha'_j|+|\beta'_j|+|\gamma'_j|=1$, $j=1,2$, $\kappa'_1,\kappa'_2\in\Omega(\beta'_1,\beta'_2)$. As above, integrating by parts with respect to $\lambda$, we obtain 
$$\left|u^{(2)}_{\alpha'_1,\alpha'_2}\right|\le C(t')^{-5/2}A_{5/2+p,5/2}(x',\xi_1,\xi_2,y'),\eqno{(5.45)}$$
where $p=|\alpha'_1|+|\alpha'_2|$. On the other hand, by (5.34), we have
$$\left|v^{(2)}_{\beta'_1,\beta'_2}\right|\le C2^{-k_1(\epsilon_1-1)}2^{-k_2(\epsilon_2-1)}|\rho_1|^{-1-\kappa_1-\epsilon_1}|\rho_2|^{-1-\kappa_2-\epsilon_2}$$ $$\times \left(|x-\xi_1|^{-\kappa_1}+|\xi_1-\xi_2|^{-\kappa_1}\right)\left(|\xi_1-\xi_2|^{-\kappa_2}+|y-\xi_2|^{-\kappa_2}\right),\eqno{(5.46)}$$
where $\kappa_j=\kappa'_j+|\beta'_1|+|\beta'_2|$, $j=1,2$. By (5.37), (5.45) and (5.46), we obtain
   $$\left|Y_{k_1,k_2}^{(4)}\right|\le C\theta^{-1} h^{-5/2}t^{-5/2}2^{-k_1(\epsilon_1-1)}2^{-k_2(\epsilon_2-1)}\sum_{p=0}^2\sum_{\kappa_1,\kappa_2\in\Omega^\sharp(p)}$$
   $$ \int_{{\bf R}^6}\int_{{\bf R}^6} |\rho_1|^{-1-\kappa_1-\epsilon_1}|\rho_2|^{-1-\kappa_2-\epsilon_2} \left(|x-\xi_1|^{-\kappa_1}+|\xi_1-\xi_2|^{-\kappa_1}\right)\left(|\xi_1-\xi_2|^{-\kappa_2}+|y-\xi_2|^{-\kappa_2}\right)$$
   $$\times A_{5/2+p,5/2}(x,\xi_1,\xi_2,y)
   \langle\xi_1\rangle^{-7/2-\epsilon'}\langle\xi_2\rangle^{-7/2-\epsilon'} d\xi_1d\xi_2,\eqno{(5.47)}$$
 where $\Omega^\sharp(p)$ denotes the set of all integers $0\le\kappa_1,\kappa_2\le 3$ such that $\kappa_1+\kappa_2\le 4-p$. Again, the integrals in the right-hand side of (5.47) are bounded by integrals satisfying the conditions of Lemma 5.6.
 \eproof
 
 If $h2^{k_1+k_2}\le 1$ we take $\theta=h2^{k_1+k_2}$ to conclude that in this case (5.36) follows from Lemma 5.7.
  If $h2^{k_1+k_2}\ge 1$ the function $Y_{k_1,k_2}$ clearly satisfies (5.38) with $\theta=1$, which again implies (5.36).

  Case 2. $|\gamma_1|+|\gamma_2|\le 1$. We will proceed as follows. If $\gamma_1=0$, $|\gamma_2|=1$, we integrate by parts once with respect to $\xi_1$, and if
$\gamma_2=0$, $|\gamma_1|=1$, we integrate by parts once with respect to $\xi_2$. When $\gamma_1=\gamma_2=0$, we integrate by parts once with respect to $\xi_1$
if $\kappa_1\le\kappa_2$, and with respect to $\xi_2$ if $\kappa_1>\kappa_2$. Then, as in the proof of (5.39) above, one can easyly see that in this case the function (5.35) satisfies the estimate
$$\left|(5.35)\right|\le C h^{-7/2}t^{-5/2}\sum_{p=0}^3\sum_{\kappa_1,\kappa_2\in\Omega^\flat(p)} \int_{{\bf R}^6}\int_{{\bf R}^6}\left(|\rho_1|^{-1}+|\rho_2|^{-1}\right)$$
   $$\times |\rho_1|^{-1-\kappa_1}|\rho_2|^{-1-\kappa_2} \left(|x-\xi_1|^{-\kappa_1}+|\xi_1-\xi_2|^{-\kappa_1}\right)\left(|\xi_1-\xi_2|^{-\kappa_2}+|y-\xi_2|^{-\kappa_2}\right)$$
   $$\times A_{5/2+p,5/2}(x,\xi_1,\xi_2,y)
   \langle\xi_1\rangle^{-7/2-\epsilon'}\langle\xi_2\rangle^{-7/2-\epsilon'} d\xi_1d\xi_2,\eqno{(5.48)}$$
 where $\Omega^\flat(p)$ denotes the set of all integers $0\le\kappa_1,\kappa_2\le 2$ such that $\kappa_1+\kappa_2\le 3-p$. It follows from (5.48) and Lemma 5.6 that in this case the function (5.35) is
 $O\left(h^{-7/2}t^{-5/2}\right)$, which is the desired result.
\eproof
  
\section*{Appendix A}

In this appendix we will sketch the proof of Lemmas 4.4, 5.3 and 5.6 following \cite{kn:EG} (see Section 6).

{\it Proof of Lemma 4.4.} It suffices to consider the integral in (4.14) in the region ${\cal O}:=\{\xi\in{\bf R}^n:|\rho(x,\xi,y)|\le\rho_0\}$, where $0<\rho_0<2$, as the bound (4.14) is trivial in $|\rho|\ge\rho_0$. Set ${\cal O}_1=\{\xi\in{\cal O}:|x-\xi|\le|y-\xi|\}$, ${\cal O}_2=\{\xi\in{\cal O}:|y-\xi|\le|x-\xi|\}$. Denote by $\xi^*$ the orthogonal projection of $\xi$ on the line $\overline{xy}$.
On ${\cal O}_1$ we introduce new coordinates $\xi=(\tau,\zeta)\in {\bf R}\times{\bf R}^{n-1}$, where
$\tau=|x-\xi^*|$ and $\zeta$ is the coordinate on the plane perpendicular to $\overline{xy}$. It is easy to see (e.g. see the proof of Theorem 3.3 of \cite{kn:EG}) that in ${\cal O}_1$ we have $|\xi-y|\sim |x-y|$, $|\xi-x|\sim\tau$, $0\le\tau\le |x-y|$,
$$|\rho(x,\xi,y)|\ge \frac{C|\zeta|}{\tau},\quad C>0,\eqno{(A.1)}$$
$$\langle\xi\rangle\sim \langle\zeta-\zeta_0\rangle+\langle\tau-\tau_0\rangle,\eqno{(A.2)}$$
where $(\tau_0,\zeta_0)$ denotes the origin in the new coordinates. We have
$$\int_{{\cal O}_1}|\rho(x,\xi,y)|^{-\ell_1}|x-\xi|^{-\ell_2}\langle\xi\rangle^{-\ell_3}d\xi$$ $$\le C\int_0^\infty\int_{|\zeta|\le\tau}\tau^{\ell_1-\ell_2}\left(\langle\zeta-\zeta_0\rangle+\langle\tau-\tau_0\rangle\right)^{-\ell_3}|\zeta|^{-\ell_1}d\zeta d\tau$$ $$=C\int_0^1 ... d\tau+C\int_1^\infty ... d\tau:=J_1+J_2.$$ 
$$J_1\le C\int_0^1\int_{|\zeta|\le\tau}\tau^{\ell_1-\ell_2}|\zeta|^{-\ell_1}d\zeta d\tau\le C\int_0^1\tau^{n-1-\ell_2}d\tau\le Const.$$
To bound $J_2$ we will consider three cases.

Case 1. $\ell_1+\ell_3<n-1$. Then we get
$$J_2\le C\int_1^\infty\int_{|\zeta|\le\tau}\tau^{\ell_1-\ell_2}|\zeta-\zeta_0|^{-\ell_3}|\zeta|^{-\ell_1}d\zeta d\tau
\le C\int_1^\infty\tau^{n-1-\ell_2-\ell_3}d\tau\le Const.$$

Case 2. $\ell_1+\ell_3=n-1$. Then $\ell_3>0$ and for every $0<\epsilon\ll 1$, we get
$$J_2\le C\int_1^\infty\int_{|\zeta|\le\tau}\tau^{\ell_1-\ell_2}|\zeta-\zeta_0|^{-\ell_3+\epsilon}|\zeta|^{-\ell_1}d\zeta d\tau
\le C\int_1^\infty\tau^{n-1-\ell_2-\ell_3+\epsilon}d\tau\le Const.$$

Case 3. $\ell_1+\ell_3>n-1$. Then, since $\ell_1\le \ell_2$, for every $0<\epsilon\ll 1$, we get
$$J_2\le C\int_1^\infty\int_{|\zeta|\le\tau}\tau^{\ell_1-\ell_2}\langle\tau-\tau_0\rangle^{n-1-\ell_1-\ell_3+\epsilon}\langle\zeta-\zeta_0\rangle^{-n+1+\ell_1-\epsilon}|\zeta|^{-\ell_1}d\zeta d\tau$$
$$\le C\int_1^\infty\tau^{\ell_1-\ell_2}\langle\tau-\tau_0\rangle^{n-1-\ell_1-\ell_3+\epsilon}d\tau$$ $$\le C\int_1^\infty\left(\tau^{n-1-\ell_2-\ell_3+\epsilon}+\langle\tau-\tau_0\rangle^{n-1-\ell_2-\ell_3+\epsilon}\right)d\tau\le Const.$$
On ${\cal O}_2$ we introduce the coordinates $\xi=(\tau,\zeta)\in {\bf R}\times{\bf R}^{n-1}$, where
$\tau=|y-\xi^*|$ and $\zeta$ is as above. In ${\cal O}_2$ we have $|\xi-x|\sim |x-y|$, $|\xi-y|\sim\tau$, $0\le\tau\le |x-y|$ as well as (A.1) and (A.2). We have as above
$$\int_{{\cal O}_2}|\rho(x,\xi,y)|^{-\ell_1}|x-\xi|^{-\ell_2}\langle\xi\rangle^{-\ell_3}d\xi$$ $$\le C|x-y|^{-\ell_2}\int_0^{|x-y|}\int_{|\zeta|\le\tau}\tau^{\ell_1}\left(\langle\zeta-\zeta_0\rangle+\langle\tau-\tau_0\rangle\right)^{-\ell_3}|\zeta|^{-\ell_1}d\zeta d\tau$$ 
$$\le C\int_0^{\infty}\int_{|\zeta|\le\tau}\tau^{\ell_1-\ell_2}\left(\langle\zeta-\zeta_0\rangle+\langle\tau-\tau_0\rangle\right)^{-\ell_3}|\zeta|^{-\ell_1}d\zeta d\tau\le Const.$$ 
\eproof

Lemma 5.3 can be derived from the following\\

\noindent
{\bf Lemma A.1} {\it Let} $ 0\le \ell'_1<n-1$, $\ell'_1\le \ell'_2<n-1$, $\ell'_2+\ell'_3>n$. 
 {\it Then we have the bound}
 $$\int_{{\bf R}^n}|\mu|^{-\ell'_1}|\xi_1-\xi_2|^{-\ell'_2}\langle\xi_2\rangle^{-\ell'_3}d\xi_2\le C\alpha(x,\xi_1,y)^{-\ell'_1},\eqno{(A.3)}$$
{\it where} $\alpha(x,\xi,y)$ {\it denotes the angle between the vectors} $\vec{\xi x}$ {\it and} $\vec{y\xi }$.\\

Indeed, in view of Lemma A.1, we obtain
$$\int_{{\bf R}^n}\int_{{\bf R}^n}|\mu|^{-\ell_1}|\xi_1-x|^{-\ell_2}|\xi_1-\xi_2|^{-\ell_3}\langle\xi_1\rangle^{-\ell_4}\langle\xi_2\rangle^{-\ell_5}d\xi_1d\xi_2$$ 
$$\le C\int_{{\bf R}^n}\alpha(x,\xi_1,y)^{-\ell_1}|\xi_1-x|^{-\ell_2}\langle\xi_1\rangle^{-\ell_4}d\xi_1.$$
In the same way as in the proof of Lemma 4.4, one can see that this integral is bounded from above by a constant.\\

{\it Proof of Lemma A.1.} We will follow the proof of Theorem 3.5 of \cite{kn:EG}. Clearly, 
 $l:=\{\xi_2:\mu=0\}=\{y+t(x-\xi_1),\,t\ge0\}$. We will only study the hardest case when
$0<\alpha\ll 1$. Denote by $K$ the cone $\{\xi: \angle(\vec{y\xi},\vec{x\xi_1})\le\alpha_0\}$, where $0<\alpha_0\ll 1$ is some constant. Clearly, in ${\bf R}^n\setminus K$ we have $|\mu|\ge\mu_0$, where $0<\mu_0<2$ is some constant. Therefore, it suffices to study the integral in $K$. On $K$ we write the variable $\xi_2$ in new coordinates $(\tau,\zeta)\in {\bf R}\times{\bf R}^{n-1}$, where
$\tau=|y-\xi_2^*|$, $\xi_2^*$ being the orthogonal projection of $\xi_2$ on the line $l$, and $\zeta$ is the coordinate on the plane perpendicular to $l$. In $K$, we have 
$$|\mu(x,\xi_1,\xi_2,y)|\ge \frac{C|\zeta|}{\tau},\quad C>0,\eqno{(A.4)}$$
as well as (A.2) with $\xi$ replaced by $\xi_2$. Let $(\tau_1,\zeta_1)$ be $\xi_1$ in the new coordinates. Clearly, $\tau_1\sim|y-\xi_1|$ and
 $|\zeta_1|\sim\alpha|y-\xi_1|$.
Thus we have
$$\int_{K}|\mu|^{-\ell'_1}|\xi_1-\xi_2|^{-\ell'_2}\langle\xi_2\rangle^{-\ell'_3}d\xi_2$$
$$\le C\int_0^{\infty}\int_{|\zeta|\le\tau}\tau^{\ell'_1}|\zeta|^{-\ell'_1}\left(|\tau-\tau_1|+|\zeta-\zeta_1|\right)^{-\ell'_2}\left(\langle\zeta-\zeta_0\rangle+
\langle\tau-\tau_0\rangle\right)^{-\ell'_3}d\zeta d\tau.\eqno{(A.5)}$$
If $\tau\ge 2|y-\xi_1|$ or $\tau\le|y-\xi_1|/2$, we have $|\tau-\tau_1|\ge\tau/C'$, $C'>1$, and in this case the right-hand side of (A.5) is bounded by
$$C\int_0^{\infty}\int_{|\zeta|\le\tau}\tau^{\ell'_1-\ell'_2}|\zeta|^{-\ell'_1}\left(\langle\zeta-\zeta_0\rangle+
\langle\tau-\tau_0\rangle\right)^{-\ell'_3}d\zeta d\tau.$$
This integral can be bounded from above by a constant in the same way as in the proof of Lemma 4.4. If $|y-\xi_1|/2\le\tau\le 2|y-\xi_1|$, we have
$|\zeta_1|\sim\alpha\tau$, and in this case the right-hand side of (A.5) is bounded by
$$C\int_{|y-\xi_1|/2}^{2|y-\xi_1|}\int_{|\zeta|\le\tau}\tau^{\ell'_1}|\zeta|^{-\ell'_1}|\zeta-\zeta_1|^{-\ell'_2}\left(\langle\zeta-\zeta_0\rangle+
\langle\tau-\tau_0\rangle\right)^{-\ell'_3}d\zeta d\tau$$
$$\le C\int_{|y-\xi_1|/2}^{2|y-\xi_1|}\int_{|\zeta|\le\tau}\tau^{\ell'_1}|\zeta_1|^{-\ell'_1}\left(|\zeta|^{-\ell'_2}+|\zeta-\zeta_1|^{-\ell'_2}\right)\left(\langle\zeta-\zeta_0\rangle+
\langle\tau-\tau_0\rangle\right)^{-\ell'_3}d\zeta d\tau$$
$$\le C\alpha^{-\ell'_1}\int_{0}^{\infty}\int_{|\zeta|\le\tau}\left(|\zeta|^{-\ell'_2}+|\zeta-\zeta_1|^{-\ell'_2}\right)\left(\langle\zeta-\zeta_0\rangle+
\langle\tau-\tau_0\rangle\right)^{-\ell'_3}d\zeta d\tau.$$
Again, this integral can be bounded as in the proof of Lemma 4.4 above.
\eproof

Lemma 5.6 is a consequence of the following\\

\noindent
{\bf Lemma A.2} {\it Let} $ 0\le \ell'_1,\ell'_2<n-1$, $0\le \ell'_3,\ell'_4<n$, $\ell'_3+\ell'_5>n$, $\ell'_4+\ell'_5>n$. {\it We also suppose that} 
$\max\{\ell'_1,\ell'_2\}\le\min\{\ell'_3,\ell'_4\}$. 
 {\it Then we have the bounds}
 $$\int_{{\bf R}^n}|\rho_1|^{-\ell'_1}|\rho_2|^{-\ell'_2}|\xi_1-x|^{-\ell'_3}\langle\xi_1\rangle^{-\ell'_5}d\xi_1\le C\alpha(x,\xi_2,y)^{-\ell'_2},\eqno{(A.6)}$$
 $$\int_{{\bf R}^n}|\rho_1|^{-\ell'_1}|\rho_2|^{-\ell'_2}|\xi_1-\xi_2|^{-\ell'_4}\langle\xi_1\rangle^{-\ell'_5}d\xi_1\le C\alpha(x,\xi_2,y)^{-\min\{\ell'_1,\ell'_2\}}.\eqno{(A.7)}$$
{\it Morever, the bounds remain true if we switch the roles of} $(\xi_1,x)$ {\it and} $(\xi_2,y)$ .\\

We will only prove (5.21) since (5.24) and (5.26) can be treated in the same way. We will first prove (5.21) under the condition (5.19). 
If $\ell_3\le\ell_4$, by (A.7), we obtain
$$\int_{{\bf R}^n}\int_{{\bf R}^n}|\rho_1|^{-\ell_1}|\rho_2|^{-\ell_2}|\xi_1-x|^{-\ell_3}|\xi_1-\xi_2|^{-\ell_4}\langle\xi_1\rangle^{-\ell_5}\langle\xi_2\rangle^{-\ell_5}d\xi_1d\xi_2$$ 
$$\le C\int_{{\bf R}^n}\alpha(x,\xi_1,y)^{-\min\{\ell_1,\ell_2\}}|\xi_1-x|^{-\ell_3}\langle\xi_1\rangle^{-\ell_5}d\xi_1.$$
If $\ell_3>\ell_4$, then $\ell_1\le\ell_3$, $\ell_2\le\ell_4$. In this case we use the inequality
$$|\xi_1-x|^{-\ell_3}|\xi_1-\xi_2|^{-\ell_4}\le |x-\xi_2|^{-\ell_4}\left(|\xi_1-x|^{-\ell_3}+|\xi_1-\xi_2|^{-\ell_3}\right)$$
together with (A.6) and (A.7) to obtain
$$\int_{{\bf R}^n}\int_{{\bf R}^n}|\rho_1|^{-\ell_1}|\rho_2|^{-\ell_2}|\xi_1-x|^{-\ell_3}|\xi_1-\xi_2|^{-\ell_4}\langle\xi_1\rangle^{-\ell_5}\langle\xi_2\rangle^{-\ell_5}d\xi_1d\xi_2$$ 
$$\le C\int_{{\bf R}^n}\alpha(x,\xi_2,y)^{-\ell_2}|\xi_2-x|^{-\ell_4}\langle\xi_2\rangle^{-\ell_5}d\xi_2.$$
In the same way as in the proof of Lemma 4.4, one can see that  these integrals are bounded from above by a constant. Suppose now (5.20) instead of (5.19). If $\min\{\ell_1,\ell_2\}\le \min\{\ell_3,\ell_4\}$, then (5.20) implies (5.19). Let now $\min\{\ell_1,\ell_2\}> \min\{\ell_3,\ell_4\}$.
In this case we use the inequality
$$|\rho_1|^{-\ell_1}|\rho_2|^{-\ell_2}\le |\rho_1|^{-\ell'_1}|\rho_2|^{-\ell'_2}+|\rho_1|^{-\ell''_1}|\rho_2|^{-\ell''_2},$$
where $\ell'_1=\min\{\ell_3,\ell_4\}$, $\ell'_2=2\ell_2-\min\{\ell_3,\ell_4\}$, $\ell''_1=2\ell_1-\min\{\ell_3,\ell_4\}$, $\ell''_2=\min\{\ell_3,\ell_4\}$. It is easy to see now that (5.20) implies that $(\ell'_1,\ell'_2)$, 
$(\ell''_1,\ell''_2)$ satisfy (5.19).\\

{\it Proof of Lemma A.2.} Again, we will follow closely the proof of Theorem 3.5 of \cite{kn:EG}. We have 
$l_1:=\{\xi_1:\rho_1=0\}=\{\xi_2+t(x-\xi_2),\,0\le t\le1\}$, $l_2:=\{\xi_1:\rho_2=0\}=\{\xi_2+t(\xi_2-y),\,t\ge0\}$. We will only study the hardest case when
$0<\alpha\ll 1$. Denote by $K_1$ (resp. $K_2$) the cone $\{\xi: \angle(\vec{\xi_2\xi},\vec{y\xi_2})\le\alpha_0\}$ (resp. $\{\xi: \angle(\vec{x\xi},\vec{x\xi_2})\le\alpha_0\}$), where $0<\alpha_0\ll 1$ is some constant. Clearly, in ${\bf R}^n\setminus K_1$ we have
$|\rho_1|\ge\rho_0$, $|\rho_2|\ge\rho_0$, while in ${\bf R}^n\setminus K_2$ we have $|\rho_1|\ge\rho_0$, where $0<\rho_0<2$ is some constant. 
We will first study the integrals in the region $K=K_1\cap K_2$. Set ${\cal O}_1^\sharp:=K\cap\{\xi: |\xi-\xi_2|\le |\xi-x|\}$, ${\cal O}_2^\sharp:=K\cap\{\xi: |\xi-x|\le |\xi-\xi_2|\}$. On ${\cal O}_1^\sharp$ we write the variable $\xi_1$ in new coordinates $(\tau,\zeta)\in {\bf R}\times{\bf R}^{n-1}$, where
$\tau=|\xi_2-\xi_1^*|$, $\xi_1^*$ being the orthogonal projection of $\xi_1$ on the line $l_2$, and $\zeta$ is the coordinate on the plane perpendicular to $l_2$. In ${\cal O}_1^\sharp$, we have $0\le\tau\le|x-\xi_2|$, $|\xi_1-\xi_2|\sim\tau$, $|\xi_1-x|\sim|\xi_2-x|$,
$$|\rho(\xi_1,\xi_2,y)|\ge \frac{C|\zeta|}{\tau},\quad C>0,\eqno{(A.8)}$$
as well as (A.2) with $\xi$ replaced by $\xi_1$. Moreover, the line $l_1$ in these coordinates can be written as $(\tau,\zeta_\tau)$ with $|\zeta_\tau|\sim\alpha\tau$, and we have
$$|\rho(x,\xi_1,\xi_2)|\ge \frac{C|\zeta-\zeta_\tau
|}{\tau},\quad C>0.\eqno{(A.9)}$$
Thus we have
$$\int_{{\cal O}_1^\sharp}|\rho_1|^{-\ell'_1}|\rho_2|^{-\ell'_2}|\xi_1-x|^{-\ell'_3}\langle\xi_1\rangle^{-\ell'_5}d\xi_1$$
$$\le C|x-\xi_2|^{-\ell'_3}\int_0^{|x-\xi_2|}\int_{|\zeta|\le\tau}\tau^{\ell'_1+\ell'_2}\left(\langle\zeta-\zeta_0\rangle+
\langle\tau-\tau_0\rangle\right)^{-\ell'_5}|\zeta|^{-\ell'_1}|\zeta-\zeta_\tau|^{-\ell'_2}d\zeta d\tau$$
$$\le C|x-\xi_2|^{-\ell'_3}\int_0^{|x-\xi_2|}\int_{|\zeta|\le\tau}\tau^{\ell'_1+\ell'_2}\left(\langle\zeta-\zeta_0\rangle+
\langle\tau-\tau_0\rangle\right)^{-\ell'_5}$$ $$\times|\zeta_\tau|^{-\min\{\ell'_1,\ell'_2\}}\left(|\zeta|^{-\max\{\ell'_1,\ell'_2\}}+|\zeta-\zeta_\tau|^{-\max\{\ell'_1,\ell'_2\}}\right)d\zeta d\tau$$
$$\le C\alpha^{-\min\{\ell'_1,\ell'_2\}}\int_0^{\infty}\int_{|\zeta|\le\tau}\tau^{\max\{\ell'_1,\ell'_2\}-\ell'_3}\left(\langle\zeta-\zeta_0\rangle+
\langle\tau-\tau_0\rangle\right)^{-\ell'_5}$$ $$\times\left(|\zeta|^{-\max\{\ell'_1,\ell'_2\}}+|\zeta-\zeta_\tau|^{-\max\{\ell'_1,\ell'_2\}}\right)d\zeta d\tau,\eqno{(A.10)}$$
$$\int_{{\cal O}_1^\sharp}|\rho_1|^{-\ell'_1}|\rho_2|^{-\ell'_2}|\xi_1-\xi_2|^{-\ell'_4}\langle\xi_1\rangle^{-\ell'_5}d\xi_1$$
$$\le C\int_0^{\infty}\int_{|\zeta|\le\tau}\tau^{\ell'_1+\ell'_2-\ell'_4}\left(\langle\zeta-\zeta_0\rangle+
\langle\tau-\tau_0\rangle\right)^{-\ell'_5}|\zeta|^{-\ell'_1}|\zeta-\zeta_\tau|^{-\ell'_2}d\zeta d\tau$$
$$\le C\alpha^{-\min\{\ell'_1,\ell'_2\}}\int_0^{\infty}\int_{|\zeta|\le\tau}\tau^{\max\{\ell'_1,\ell'_2\}-\ell'_4}\left(\langle\zeta-\zeta_0\rangle+
\langle\tau-\tau_0\rangle\right)^{-\ell'_5}$$ $$\times\left(|\zeta|^{-\max\{\ell'_1,\ell'_2\}}+|\zeta-\zeta_\tau|^{-\max\{\ell'_1,\ell'_2\}}\right)d\zeta d\tau.\eqno{(A.11)}$$
The integral in the right-hand side of (A.10) and (A.11) can be bounded from above by a constant in the same way as in the proof of Lemma 4.4. The integrals over ${\cal O}_2^\sharp$ can be studied similarly with the main difference that in this case we take $\tau=|x-\xi_1^{**}|$, $\xi_1^{**}$ being the orthogonal projection of $\xi_1$ on the line $l_1$, and $\zeta$ is the coordinate on the plane perpendicular to $l_1$. In this case we have $0\le\tau\le|x-\xi_2|$, $|\xi_1-x|\sim\tau$, $|\xi_1-\xi_2|\sim|\xi_2-x|$,
$$|\rho(x,\xi_1,\xi_2)|\ge \frac{C|\zeta|}{\tau},\quad C>0.\eqno{(A.12)}$$
And if the line $l_1$ in these coordinates is written as $(\tau,\zeta_\tau)$, we also have $|\zeta_\tau|\sim\alpha|\xi_2-x|$ and 
$$|\rho(\xi_1,\xi_2,y)|\ge \frac{C|\zeta-\zeta_\tau
|}{\tau},\quad C>0.\eqno{(A.13)}$$
It remains to bound our integrals in the region $K_1\setminus K$. We write the variable $\xi_1$ in new coordinates $(\tau,\zeta)\in {\bf R}\times{\bf R}^{n-1}$, where
$\tau=|\xi_2-\xi_1^*|$, $\xi_1^*$ being the orthogonal projection of $\xi_1$ on the line $l_2$, and $\zeta$ is the coordinate on the plane perpendicular to $l_2$. In $K_1\setminus K$, we have $\tau\ge|x-\xi_2|$, $|\xi_1-\xi_2|\sim\tau$, and 
$$|\rho(x,\xi_1,\xi_2)|\ge \rho_0>0,\eqno{(A.14)}$$
$$|\rho(\xi_1,\xi_2,y)|\ge \frac{C|\zeta|}{\tau},\quad C>0,\eqno{(A.15)}$$
as well as (A.2) with $\xi$ replaced by $\xi_1$. Let $(\tau_x,\zeta_x)$ be $x$ in the new coordinates. Clearly, $\tau_x\sim|x-\xi_2|$ and
 $|\zeta_x|\sim\alpha|x-\xi_2|$.
Thus we have
$$\int_{K_1\setminus K}|\rho_2|^{-\ell'_2}|\xi_1-\xi_2|^{-\ell'_4}\langle\xi_1\rangle^{-\ell'_5}d\xi_1$$
$$\le C\int_0^{\infty}\int_{|\zeta|\le\tau}\tau^{\ell'_2-\ell'_4}|\zeta|^{-\ell'_2}\left(\langle\zeta-\zeta_0\rangle+
\langle\tau-\tau_0\rangle\right)^{-\ell'_5}d\zeta d\tau\le Const,\eqno{(A.16)}$$
$$\int_{K_1\setminus K}|\rho_2|^{-\ell'_2}|\xi_1-x|^{-\ell'_3}\langle\xi_1\rangle^{-\ell'_5}d\xi_1$$
$$\le C\int_{|x-\xi_2|}^{\infty}\int_{|\zeta|\le\tau}\tau^{\ell'_2}|\zeta|^{-\ell'_2}\left(|\tau-\tau_x|+|\zeta-\zeta_x|\right)^{-\ell'_3}\left(\langle\zeta-\zeta_0\rangle+
\langle\tau-\tau_0\rangle\right)^{-\ell'_5}d\zeta d\tau$$
$$=C\left(\int_{2|x-\xi_2|}^{\infty}+\int_{|x-\xi_2|}^{2|x-\xi_2|}\right)\int_{|\zeta|\le\tau}\tau^{\ell'_2}|\zeta|^{-\ell'_2}\left(|\tau-\tau_x|+|\zeta-\zeta_x|\right)^{-\ell'_3}\left(\langle\zeta-\zeta_0\rangle+
\langle\tau-\tau_0\rangle\right)^{-\ell'_5}d\zeta d\tau$$ $$:=I_1+I_2.\eqno{(A.17)}$$
If $\tau\ge 2|x-\xi_2|$, we have $|\tau-\tau_x|\ge\tau/C'$, $C'>1$, and the integral $I_1$ is bounded by a constant as in (A.16).
 If $|x-\xi_2|\le\tau\le 2|x-\xi_2|$, we have $|\zeta_x|\sim\alpha\tau$. To bound $I_2$ we will consider two cases:
 
 Case 1. $0\le\ell'_3<n-1$. Then 
$$I_2\le  C\int_{|x-\xi_2|}^{2|x-\xi_2|}\int_{|\zeta|\le\tau}\tau^{\ell'_2}|\zeta_x|^{-\ell'_2}\left(|\zeta|^{-\ell'_3}+|\zeta-\zeta_x|^{-\ell'_3}\right)\left(\langle\zeta-\zeta_0\rangle+
\langle\tau-\tau_0\rangle\right)^{-\ell'_5}d\zeta d\tau$$
$$\le C\alpha^{-\ell'_2}\int_{0}^{\infty}\int_{|\zeta|\le\tau}\left(|\zeta|^{-\ell'_3}+|\zeta-\zeta_x|^{-\ell'_3}
\right)\langle\zeta-\zeta_0\rangle^{-n+1+\ell'_3}
\langle\tau-\tau_0\rangle^{-\ell'_5-\ell'_3+n-1}d\zeta d\tau$$
$$\le C\alpha^{-\ell'_2}\int_{0}^{\infty}
\langle\tau-\tau_0\rangle^{-\ell'_5-\ell'_3+n-1}d\tau\le C\alpha^{-\ell'_2}.$$

Case 2. $n-1\le\ell'_3<n$. Then
$$I_2\le C\int_{|x-\xi_2|}^{2|x-\xi_2|}\int_{|\zeta|\le\tau}\tau^{\ell'_2}|\zeta|^{-\ell'_2}|\zeta-\zeta_x|^{-n+1-\epsilon}|\tau-\tau_x|^{-\ell'_3+n-1-\epsilon}\left(\langle\zeta-\zeta_0\rangle+
\langle\tau-\tau_0\rangle\right)^{-\ell'_5}d\zeta d\tau$$
$$\le  C\int_{|x-\xi_2|}^{2|x-\xi_2|}\int_{|\zeta|\le\tau}\tau^{\ell'_2}|\zeta_x|^{-\ell'_2}\left(|\zeta|^{-n+1-\epsilon}+|\zeta-\zeta_x|^{-n+1-\epsilon}\right)$$ $$\times|\tau-\tau_x|^{-\ell'_3+n-1-\epsilon}\left(\langle\zeta-\zeta_0\rangle+
\langle\tau-\tau_0\rangle\right)^{-\ell'_5}d\zeta d\tau$$
$$\le C\alpha^{-\ell'_2}\int_{0}^{\infty}\int_{|\zeta|\le\tau}\left(|\zeta|^{-n+1-\epsilon}+|\zeta-\zeta_x|^{-n+1-\epsilon}\right)\langle\zeta-\zeta_0\rangle^{-2\epsilon}$$ $$\times|\tau-\tau_x|^{-\ell'_3+n-1-\epsilon}
\langle\tau-\tau_0\rangle^{-\ell'_5+2\epsilon}d\zeta d\tau$$
$$\le C\alpha^{-\ell'_2}\int_{0}^{\infty}|\tau-\tau_x|^{-\ell'_3+n-1-\epsilon}
\langle\tau-\tau_0\rangle^{-\ell'_5+2\epsilon} d\tau\le C\alpha^{-\ell'_2}.$$
\eproof

\section*{Appendix B}

To prove (2.1) (with $\epsilon=0$) in all dimensions $n\ge 4$, one is led to bound oscilatory integrals of the kind
$$I_k(h,m)=\int_{{\bf R}^n}...\int_{{\bf R}^n}e^{i\psi/h}a^mV(\xi_1)...V(\xi_k)d\xi_1...d\xi_k,$$
where
$$\psi=|x-\xi_1|+|\xi_1-\xi_2|+...+|\xi_{k-1}-\xi_k|+|\xi_k-y|,$$
$$a=\frac{|x-\xi_1|+|\xi_1-\xi_2|+...+|\xi_{k-1}-\xi_k|+|\xi_k-y|}{|x-\xi_1||\xi_1-\xi_2|...|\xi_{k-1}-\xi_k||\xi_k-y|},$$
$x,y\in{\bf R}^n$, $k\ge 1$, $1\le m <n$ and $0<h\ll 1$. The key point is the following estimate which seems hard to prove for $k\ge 2$ and all
$n\ge 4$.\\

\noindent
{\bf Conjecture.} {\it If $V$ satisfies (1.6), then} 
$$\left|I_k\left(h,\frac{n-1}{2}\right)\right|\le Ch^{k(n-3)/2},$$
{\it with a constant $C>0$ independent of $x,y$ and $h$.}\\

Note that when $n=5$ this conjecture is actually proved in \cite{kn:EG}.\\

{\bf Acknowledgements.}  F.C. has been partially supported by the CNPq-Brazil, while G.V. has been partially supported by the project NONAa, ANR-08-BLAN-0228. Both authors have also been partially supported by the agreement Brazil-France in Mathematics - Proc. 69.0014/01-5.

F. Cardoso

Universidade Federal de Pernambuco, 

Departamento de Matem\'atica, 

CEP. 50540-740 Recife-Pe, Brazil,

e-mail: fernando@dmat.ufpe.br\\

G. Vodev

Universit\'e de Nantes,

 D\'epartement de Math\'ematiques, UMR 6629 du CNRS,
 
 2, rue de la Houssini\`ere, BP 92208, 
 
 44332 Nantes Cedex 03, France,
 
 e-mail: vodev@math.univ-nantes.fr


\begin{thebibliography}
\frenchspacing \baselineskip=12 pt plus 1pt minus 1pt

\bibitem{kn:B} {\sc M. Beals}, {\em Optimal $L^\infty$ decay estimates for solutions to the wave equation with a potential}, Comm. PDE {\bf 19} (1994), 1319-1369.

\bibitem{kn:CV} {\sc F. Cardoso and G. Vodev},
{\em Semi-classical dispersive estimates for the wave and Schr\"odinger
equations with a potential in dimensions $n\ge 4$}, Cubo Math. J. {\bf 10} (2008), 1-14.

\bibitem{kn:CCV1} {\sc F. Cardoso, C. Cuevas and G. Vodev}, {\em Dispersive estimates of solutions to the wave
equation with a potential in dimensions two and three}, Serdica Math. J. {\bf 31} (2005), 263-278.

\bibitem{kn:CCV2} {\sc F. Cardoso, C. Cuevas and G. Vodev},
{\em Dispersive estimates for the Schr\"odinger
equation in dimensions four and five}, Asymptot. Anal. {\bf 62} (2009), 125-145.

\bibitem{kn:CCV3} {\sc F. Cardoso, C. Cuevas and G. Vodev},
{\em Dispersive estimates for the Schr\"odinger
equation with potentials of critical regularity}, Cubo Math. J. {\bf 11} (2009), 57-70.

\bibitem{kn:DP} {\sc P. D'ancona and V. Pierfelice}, {\em On the wave equation with a large rough potential}, J. Funct. Analysis {\bf 227} (2005), 30-77.

\bibitem{kn:EG} {\sc M. B. Erdogan and W. R. Green}, {\em Dispersive estimates for the Schr\"odinger equation for $C^{\frac{n-3}{2}}$ potentials in odd dimensions}, IMRN {\bf 2010} (2010), No. 13, 2532-2565.

\bibitem{kn:GV} {\sc V. Georgiev and N. Visciglia}, {\em Decay estimates for the wave equation with potential}, Comm. PDE {\bf 28} (2003), 1325-1369.

\bibitem{kn:GoV} {\sc M. Goldberg and M. Visan},
{\em A counterexample to dispersive estimates for Schr\"odinger
operators in higher dimensions}, Commun. Math. Phys. {\bf 266}
(2006), 211-238.

\bibitem{kn:M} {\sc S. Moulin},
{\em Low-frequency dispersive estimates for the wave equation in higher dimensions}, Asymptot. 
Anal. {\bf 60} (2008), 15-27.

\bibitem{kn:M2} {\sc S. Moulin},
{\em High frequency dispersive estimates in dimension two}, Ann. H. Poincar\'e {\bf 10} (2009),
415-428.

\bibitem{kn:V} {\sc G. Vodev},
{\em Dispersive estimates of solutions to the wave equation with a potential 
in dimensions $n\ge 4$}, Comm. PDE {\bf 31} (2006), 1709-1733.


\end{thebibliography}
\end{document}